\documentclass{article}

\usepackage{arxiv}
\usepackage[utf8]{inputenc}
\usepackage[T1]{fontenc}
\usepackage{hyperref}
\usepackage{url}
\usepackage{booktabs}
\usepackage{amsfonts}
\usepackage{nicefrac}
\usepackage{microtype}
\usepackage{lipsum}
\usepackage{amsmath,amssymb}
\usepackage{amsthm}
\usepackage{graphicx}
\usepackage{listings}
\usepackage{color}
\usepackage{fancyhdr}
\usepackage[framemethod=tikz]{mdframed}
\usepackage{subfig}
\usepackage[nottoc,notlot,notlof]{tocbibind}
\usepackage{csvsimple}
\usepackage{longtable}
\usepackage{colortbl}
\usepackage{float}
\usepackage[linesnumbered,commentsnumbered,ruled]{algorithm2e}
\usepackage{lscape}
\usepackage{multirow}
\usepackage{hhline}
\usepackage[backend=biber,style=numeric,sorting=anyt,url=false]{biblatex}
\usepackage{arydshln}
\usepackage{nicematrix}
\usepackage{todonotes}
\usepackage[draft]{changes}
\usepackage{tikz}

\usetikzlibrary{arrows.meta}

\definechangesauthor[name={Roland Masson}, color=orange]{RM}
\definechangesauthor[name={Thomas Jacumin}, color=violet]{TJ}
\definechangesauthor[name={Laurent Jeannin}, color=red]{LJ}
\definechangesauthor[name={Arthur Dartois}, color=green]{AD}

\title{A Coats-Based Formulation for Multiphase Reactive Transport Accounting for Aqueous Phase Disappearance}

\renewcommand{\arraystretch}{1.2}
\usetikzlibrary{decorations}
\NiceMatrixOptions
{
	custom-line =
	{
		letter = ; ,
		command = \hdashedline ,
		ccommand = \cdashedline ,
		tikz = { dash pattern = on 4 pt off 2pt } ,
		total-width = \pgflinewidth
	}
}

\author{
	Thomas \textsc{Jacumin} \\
	Université Côte d’Azur, Inria, CNRS, LJAD, Team Galets, Nice, France \\
	\texttt{thomas.jacumin@inria.fr} \\
	\and
	\textbf{Roland \textsc{Masson}} \\ 
	Université Côte d’Azur, Inria, CNRS, LJAD, Team Galets, Nice, France \\
	\texttt{roland.masson@inria.fr} \\
	\and
	\textbf{Laurent \textsc{Jeannin}} \\ 
	STORENGY, Direction Expertise Technique, La Garenne-Colombes, France \\
	\texttt{laurent.jeannin@storengy.com} \\
	\and
	\textbf{Arthur \textsc{Dartois}} \\ 
	STORENGY, Direction Expertise Technique, La Garenne-Colombes, France \\
	\texttt{arthur.dartois@storengy.com} \\
}

\def\R{\mathbb{R}}
\DeclareMathOperator{\divergence}{div}
\DeclareMathOperator{\rank}{rank}
\DeclareMathOperator{\image}{Im}

\def\E{\mathcal{E}}
\def\C{\mathcal{C}}
\def\P{\mathcal{P}}
\def\Q{\mathcal{Q}}
\def\Rea{\mathcal{R}}
\def\l{\ell}
\def\x{\mathbf{x}}
\def\u{\mathbf{u}}
\def\1{\mathbf{1}}
\def\cells{\mathcal{M}}
\def\faces{\mathcal{F}}

\def\cod{{\text{CO}_2}}
\def\codg{{\text{CO}_2^{g}}}
\def\codl{{\text{CO}_2^{\l}}}
\def\hp{{\text{H}^{+}}}
\def\ohm{{\text{OH}^{-}}}
\def\hdo{{\text{H}_2\text{O}}}
\def\hdog{{\text{H}_2\text{O}^{g}}}
\def\hdol{{\text{H}_2\text{O}^{\l}}}
\def\hcotm{{\text{HCO}_3^{-}}}
\def\cacot{{\text{CaCO}_3}}
\def\cadp{{\text{Ca}^{2+}}}

\def\cotdm{{\text{CO}_3^{2-}}}
\def\hp{{\text{H}^+}}
\def\ohm{{\text{OH}^-}}
\def\soqdm{\text{SO}_4^{2-}}
\def\hdl{\text{H}_2^\l}
\def\hd{\text{H}_2}
\def\hsm{\text{HS}^-}
\def\hdsl{\text{H}_2\text{S}^\l}
\def\hdsg{\text{H}_2\text{S}^g}
\def\hdg{\text{H}_2^g}
\def\nap{\text{Na}^+}

\begin{document}
	
	\maketitle
	
	\begin{abstract}
		As part of the energy transition, underground hydrogen ($\hd$) storage in porous media
		has emerged as an alternative solution for large-scale energy storage to salt caverns. Assessing the feasibility of such systems and ensuring their operation requires accurate modeling of the geochemical processes that may alter the composition of the stored gas. This results in a highly nonlinear coupled system of partial differential equations, ordinary differential equations, and algebraic constraints, calling for robust and specialized numerical strategies.
		
		In this work, we adopt a Globally Implicit Approach (GIA), as sequential solution algorithms are known to introduce splitting errors due to the decoupling of flow and geochemical processes \cite{trenty_sequential_2006}. In addition, most splitting-based methods implicitly assume the persistent presence of the liquid phase, an assumption that may be violated during gas injection or withdrawal owing to aqueous phase vaporization and mineral precipitation in the vicinity of  wells.
		%
		We extend Coats' formulation to multiphase reactive transport
		with a rigorous and fully implicit management of aqueous phase disappearance and reappearance. The proposed approach introduces new molar concentration unknowns for basis species that track the re-appearance of the aqueous phase and guarantee the nonsingularity of the Jacobian matrix. Additionally, the appearance of the liquid phase is managed using a multiphase flash calculation algorithm \cite{kala_element_2020}.
		
		The proposed algorithm is validated against the $\cod$ storage benchmark problems \cite{ahusborde_benchmark_2024,de_hoop_benchmark_2024}. Finally, the approach is applied to an underground $\hd$ storage in an aquifer under cyclic injection and withdrawal conditions.
	\end{abstract}
	\keywords{Hydrogen storage, Multiphase reactive transport, Globally Implicit Approach, Phase transition}
	        
	\section*{Introduction}
	In the context of the energy transition, underground gas storage operators are investigating and developing solutions for underground hydrogen ($\hd$) storage in porous media. As part of the development of the hydrogen sector in Europe, pure $\hd$ produced from renewable electricity via electrolysis can contribute to the decarbonization of electricity generation and certain industrial processes. When stored underground at a large scale, hydrogen serves both as a buffer for industrial uses and for electricity networks by helping to balance production and demand.
	
	The objective of this study is to investigate the behavior of underground $\hd$ storage in aquifers. One of the challenges is to account for chemical processes likely to modify the gas composition within the storage. For example, methanation or sulfate reduction reactions, linked to microbial activity, can enrich the gas in methane or hydrogen sulfide. This work aims to numerically model two-phase water-gas flow in porous media while accounting for geochemical or biogeochemical reactions in the medium.	
	
	The resolution of these systems presents numerous difficulties related to the large number of chemical species, the stiffness of the system induced by chemical reactions, and the degeneracies caused by phase transitions.
	Most existing formulations rely on splitting algorithms that solve the coupled problem by decomposing it into two well-established submodels, coupled sequentially or iteratively. The first submodel typically consists of a compositional multiphase flow model with equilibrium phase-change reactions (the RS step). The second is a reactive transport model solved for prescribed phase velocities and volume fractions (the RT step), usually involving the mobile aqueous phase and mineral phases. Several splitting strategies have been proposed in the literature. For instance, \cite{trenty_sequential_2006} partitions the reactive system between RS and RT reactions, whereas \cite{ahusborde_numerical_2015,ahusborde_improvement_2018} partition the chemical species between RS and RT subsystems. An alternative approach is proposed in \cite{sin_integrating_2017}, where the RT step accounts for all reactions while phase velocities and volume fractions are provided by the RS flow model. However, these methods suffer from splitting errors arising from the loss of molar and volumetric conservation, which may severely restrict the allowable time step size. In addition, they generally assume the persistent presence of the aqueous phase, thereby limiting their ability to model drying phenomena occurring, for instance, near gas injection wells.
	
	To overcome these limitations, this work focuses on fully coupled approaches, also referred to as \textit{Globally Implicit Approaches} (GIA). Several GIA algorithms have been proposed in the literature. For instance the formulation proposed in \cite{nghiem_modeling_2004} relies on the species number of moles per unit reference volume as main unknowns in addition to the pressure combined with an elimination of the reaction rates using total primary species conservation equations. The work of \cite{fan_fully-coupled_2012} extends the Coats' formulation widely used in reservoir simulation \cite{coats-89-imp,EGHM2012} to reactive transport problems. This approach relies on so-called natural variables (pressure, saturations, and species mole fractions in each phase) and reduces the number of flash calculations at the expense of variable switching during the Newton iterations \cite{Voskov2011}.
	In \cite{farshidi_compositional_2016}, an alternative to natural-variable formulations is proposed to improve robustness in situations involving aqueous-phase disappearance near injection wells. The method is based on total component mole fractions in the fluid phases.
	The formulation introduced in \cite{kala_element_2020} employs total mole fractions of primary elements together with a multiphase reactive flash calculation. To reduce computational cost, the flash outputs are approximated using multilinear interpolation within the \textit{Operator-Based Linearization} framework.
	The approaches developed in \cite{ahusborde_numerical_2021,ahusborde_numerical_2024} are based on a Morel-type formulation using conservation equations for total primary species, thereby eliminating reaction rates from the governing equations. Secondary-species mole fractions are removed at the nonlinear level through mass-action laws and phase-equilibrium relations. Finally, let us mention \cite{knodel_global_2022} which extends to multiphase systems the reduction strategy introduced by \cite{krautle_reduction_2007}, yielding a stronger decoupling of the governing equations than standard formulations such as the Morel formalism.	
	
	The main objective of this work is to extend natural-variable formulations to reactive systems  while addressing the challenges associated with aqueous-phase disappearance and multiphase chemical equilibria. The proposed approach introduces additional molar concentration unknowns for primary species in absent phases. These variables are designed to preserve the nonsingularity of the Jacobian matrix and to enable accurate tracking of aqueous-phase reappearance.
	The formulation is first validated on the two-phase reactive transport benchmark proposed in \cite{ahusborde_benchmark_2024,de_hoop_benchmark_2024}, which exhibits aqueous-phase disappearance in the vicinity of a $\cod$ injection boundary. It is then applied to the simulation of hydrogen injection and production cycles into a porous storage reservoir.
	
	The remainder of this paper is organized as follows. Section~\ref{sec:model} presents the reactive transport model considered throughout the paper, based on Morel's formulation of chemical reactions. Section~\ref{sec:gia} introduces the Coats-type formulation combined with a cell-centered finite-volume discretization in space and a fully implicit Euler scheme in time. Particular attention is devoted to the treatment of absent phases through the introduction of additional molar concentration variables, as well as to the handling of phase appearance and disappearance. Finally, Section~\ref{sec:num} presents numerical test cases assessing the validity and efficiency of the proposed formulation on both a two-dimensional $\cod$ benchmark and a hydrogen storage scenario in a porous reservoir under cyclic injection and withdrawal conditions.
	
	\section{Reactive transport model}
	\label{sec:model}
	
	In this section, we present the multiphase reactive transport model, which leads to a system of partial differential equations coupled with linear and nonlinear algebraic constraints. Note that a table containing our notations for the main physical quantities used throughout the article is given in Appendix~\ref{sec:notations}. Let $\Omega$ be an open bounded domain of $\R^d$, $d=\{1,2,3\}$ and $(0, T_\text{end})$ with $T_\text{end} > 0$ be the simulation time interval. The chemical system $\mathcal{S} = \{\E, \Rea, \P\}$ is described by a collection of three sets:
	\begin{itemize}
		\item the set $\E$ of chemical species of cardinality denoted by $N$,
		\item the set  $\P = \P^\text{f}\cup \P^\text{s}$ of phases with $\P^\text{f}$ the set of fluid phases and $\P^\text{s}$ the set of solid phases, 
		\item the set $\Rea$ of the $N-M$ chemical reactions. 	
	\end{itemize}
	
	To link a species $e\in \E$ to its phase $\beta\in \P$, we introduce the following map:
	$$ \alpha : e\in\E \mapsto \alpha(e)\in\P. $$
	In this way, $\alpha^{-1}\big(\{\beta\}\big)$ represents the set of all species in the phase $\beta$. Note that a solid phase is assumed to be associated to a unique species in $\E$. In the following, the notations
	$$ \E^\text{f} = \alpha^{-1}\big( \P^\text{f} \big) \quad \mbox{and} \quad \E^\text{s} = \alpha^{-1}\big( \P^\text{s} \big), $$
	will be used to denote respectively the sets of fluid and solid species. 
	
	\subsection{Morel's formalism and total molar conservation equations}
	
	Following Morel's formalism \cite{morel_numerical_1972}, we aim to cancel out the reaction rates by using well-chosen linear combinations of the molar conservation equations defining the total molar conservation equations. For that, we introduce the set of primary species
	$$ \E^\text{pr} = \{ e^\text{pr}_1, \hdots, e^\text{pr}_M \} \subset\E, $$
	and the set of secondary species
	$$ \E^\text{sec} = \{ e^\text{sec}_1,\hdots, e^\text{sec}_{N-M} \}\subset\E, $$
	such that $(\E^\text{pr},\E^\text{sec})$ is a partition of $\E$.
	Then, every species $e$ in $\E$ is assumed to be expressed as a linear combination of the primary species as follows: 
	$$ e = \sum_{e^\text{pr}\in\E^\text{pr}} E_{e^\text{pr},e}\, e^\text{pr}, $$
	where
	$$ E := [E_{e^\text{pr},e}]_{e^\text{pr}\in\E^\text{pr},e\in\E}\in\R^{\E^\text{pr}\times\E}, $$
	is the formula matrix whose column associated to $e \in \E$ defines the composition of the species $e$ in terms of the primary species $e^\text{pr} \in \E^{\text{pr}}$, and each secondary species $e^\text{sec} \in \E^{\text{sec}}$ is associated to the chemical reaction
	$$ e^\text{sec} \leftrightarrows \sum_{e^\text{pr}\in\E^\text{pr}} E_{e^\text{pr},e^\text{sec}}\, e^\text{pr}.  $$
	As a result, we shall use in the sequel the notation $e^\text{sec}\in\Rea$ to identify the secondary species with its reaction in $\Rea$. Note that if $e$ is a primary species, then
	$$ E_{e^\text{pr},e} = \begin{cases}
		1 & \text{if}\ e^\text{pr} = e, \\
		0 & \text{otherwise}.
	\end{cases} $$
	Let us notice that, by definition, $E$ has full rank equal to $M$. 
	The elimination of the reaction rates classically leads to the following system of molar conservation of total primary species: 
	\begin{equation}
		\sum_{{e}\in \E} E_{e^{\text{pr}},e} \partial_t n_{e} + \sum_{{e}\in \E^\text{f}} E_{e^{\text{pr}},e} \divergence\Big(\zeta^{\alpha(e)} C_{e} \u^{\alpha(e)}\Big) = 0, \quad e^{\text{pr}} \,\in\, \mathcal{E}^{\text{pr}}. 
		\label{eq:esp-base}
	\end{equation}
	where $n_e$ is the volumetric molar concentration of $e$ [$\mathrm{mol\cdot m^{-3}}$], considered as an independent variable for $e\in \E^\text{s}$ and equal to
	$$ n_e = \phi \zeta^{\alpha(e)} S^{\alpha(e)} C_e, $$
	for fluid species $e\in \E^{\text{f}}$,
	with $\phi$ the porosity or pore volume fraction [$-$],
	$\zeta^\beta$ the phase molar density  [$\mathrm{mol\cdot m^{-3}}$],
	$S^\beta$ the phase saturation [$-$],
	and $C_e$ the mole fraction of $e$ in its phase [$-$]. The velocity of the fluid phase $\beta$ [$\mathrm{m\cdot s^{-1}}$] is denoted by 
	$\mathbf{u}^\beta$ and defined by Darcy's law, in which we neglect the diffusion effect 
	$$ \mathbf{u}^\beta = - \mathbb{K}\, \eta^\beta \left(\nabla P^\beta + \rho^\beta\mathbf{g} \right), $$
	with $\mathbb{K}$ the permeability tensor [$\mathrm{m^2}$], $\eta^\beta := \frac{k_\text{r}^\beta}{\mu^\beta}$ the phase mobility [$\mathrm{Pa^{-1} \cdot s^{-1}}$], $k_\text{r}^\beta$ the phase relative permeability [$-$], $\mu^\beta$ the phase viscosity [$\mathrm{Pa \cdot s}$], $P^\beta$ the phase pressure [$\mathrm{Pa}$], $\rho^\beta$ the phase mass density [$\mathrm{kg\cdot m^{-3}}$], $\mathbf{g}$ the gravitational acceleration vector [$\mathrm{m\cdot s^{-2}}$]. In the sequel, for simplicity,  the capillary pressure is neglected i.e.
	$$ P^\beta = P,\quad \forall\beta\in\P^\mathrm{f}. $$
	
	\subsection{Chemical reactions, volume conservations and other closure laws}
	
	Each chemical reaction or equivalently each secondary species $e\in \E^\text{sec}$ is associated with one closure equation which differs for equilibrium and kinetic reactions. 
	Phase transition reactions are all assumed at equilibrium which is enforced through fugacity equality, whereas other chemical equilibrium such as those for homogeneous speciation or heterogeneous precipitation/dissolution is described using the law of mass action. Each kinetic reaction for ``kinetic species'' $e$ is associated to a conservation equation with specified kinetic reaction rate $T_e$:
	\begin{align*}
		\partial_t n_e + \divergence\Big(\zeta^{\alpha(e)} C_{e} \u^{\alpha(e)}\Big) &= T_e\ \mbox{ if e is a fluid kinetic species}, \\
		\partial_t n_e &= T_e\  \mbox{ if e is a solid kinetic species}. 
	\end{align*}
	
	The total volume conservation of fluid and minerals is written as follows 
	$$
	\underbrace{\phi}_{\text{fluid}} + \underbrace{\phi^\text{s} + 1 - \phi_\text{T}}_{\text{solid}} = 1
	\quad \Leftrightarrow \quad 
	\phi + \phi^\text{s} = \phi_\text{T},
	$$
	where $\phi^\text{s}$ denotes the volume fraction of reactive rock [$-$], and $1 - \phi_\text{T}$ represents the volume fraction of non-reactive rock [$-$], which therefore remains invariant under chemical reactions.
	
	Finally, the system is closed by the fluid pore volume conservation:  
	\begin{equation} \label{eq:saturated-pores}
		\sum_{\beta\in\P^\text{f}} S^\beta = 1, 
	\end{equation}
	the sum to one of each fluid phase  molar fractions: 
	\begin{equation} \label{eq:close-comp-fluid}
		\sum_{e \in \alpha^{-1}(\{\beta\})} C_e = 1, \quad \beta\in\P^\text{f},
	\end{equation}
	and the definition of the reactive solid porosity \cite{de_hoop_benchmark_2024}
	\begin{equation} \label{eq:close-porosity}
		\phi^\text{s} = \sum_{m \in \P^\text{s}} \frac{n_m}{\zeta^m}.
	\end{equation}
	
	\subsection{Summary of the set of equations and unknowns}
	
	To conclude this section, let us summarize the system composed of $N + \# \P^\text{f}+ 2$ equations: 
	\begin{itemize}
		\item  the conservation equations for total primary species \eqref{eq:esp-base}, 
		\item  one algebraic or differential equation for each equilibrium or kinetic chemical reaction, respectively,
		\item the closure relations \eqref{eq:saturated-pores}, \eqref{eq:close-comp-fluid}, and \eqref{eq:close-porosity}.
	\end{itemize}
	These equations are associated with the $N + \# \P^\text{f}+ 2$ so-called natural variables: 
	\begin{itemize}
		\item  pressure $P$,  
		\item  saturations $S^\beta, \quad \beta \in \P^\text{f}$,
		\item solid porosity $\phi^\text{s}$,  
		\item  molar fractions $C_e, \quad e\in \E^\text{f}$,
		\item molar concentrations $n_e, \quad e\in \E^\text{s}$.  
	\end{itemize}
	This counting argument, however, assumes that all phases are present throughout the space-time domain. Accounting for the actual set of present phases is a key aspect of natural-variable formulations, since mole fractions are defined only for existing phases. This requires a careful adaptation of both the set of unknowns and of governing equations, as detailed in the next section.

	\section{Coats' formulation of reactive transport model with phase appearance and disappearance}\label{sec:gia}
	
	In this section the Coats' formulation is extended to the reactive transport model taking into account phase appearance and disappearance. To fix ideas, we consider a cell-centered finite volume discretization in space \cite{eymard_finite_2000} combined with an Euler implicit time integration. 
	
	\subsection{Space and time discretizations}
	
	Let $\cells_h$ denote a set of polytopal cells $K$ defining a conforming partition of the domain $\Omega$ and denote by $|K|$ the volume of the cell. The set of interior faces will be denoted by $\faces_h$ and the set of faces of a given cell $K$ by $\faces_K$. It is assumed that each interior face $\sigma$ is shared by exactly two cells $K$ and $L$ using the notation $\sigma=K|L$. Homogeneous Neumann boundary conditions are assumed for simplicity at the boundary of $\Omega$ in order to omit the boundary fluxes. We denote by $V_{KL}^\beta$ a discretization  of the fluid phase $\beta$ Darcy flux (excluding phase mobility) at face $\sigma = K|L$
	$$ \int_{\sigma} -\mathbb{K} \big( \nabla P + \rho^\beta \mathbf{g} \big) \cdot \mathbf{n}_{KL}, $$
	assumed to be conservative such that $V_{KL}^\beta + V_{LK}^\beta = 0$.  The mobility terms $\zeta^\beta C_e \eta^\beta$ appearing in the fluxes at face $\sigma=K|L$ are discretized using the upwind value $(\zeta^\beta C_e \eta^\beta)_{K^\beta_\sigma}$ defined by a phase-based upwind scheme \cite{aziz-settari-79,peaceman-77} with 
	$$
	K_\sigma^\beta = \left\{ \begin{array}{ll}
		K & \text{if}\ V_{KL}^\beta \geq 0, \\
		L & \text{if}\ V_{KL}^\beta < 0.
	\end{array}\right.
	$$
	The time discretization is defined by the sequence of $m+1$ discrete times
	$$
	0 = t^0 < t^1 < \cdots < t^{n-1} < t^n < \cdots < t^m=T_\text{end},
	$$ 
	and we denote by $\Delta t^n = t^n-t^{n-1}$ the time step $n$. Note that in the sequel, the current time step $n$ superscript will be most of time omitted while the previous time step $n-1$ is always kept. 
	
	\subsection{Coats' formulation}	
	
	To account for the potential disappearance and appearance of phases as a result of chemical reactions, let us introduce as additional variable the subset of present phases 
	$$
	\Q = (\Q^\text{f}, \Q^\text{s}),
	$$
	where $\Q^{\text{f}}\subset \P^\text{f}$ (resp. $\Q^{\text{s}}\subset \P^\text{s}$) represents the subset of present fluid  (resp. solid) phases. The complementary set of absent phases is denoted by 
	$$
	\Q^{c} =  \P \setminus \Q. 
	$$
	As usual for the natural variable or Coats' formulation, the set of unknowns in cell $K$ and at time step $n$ depends on the set of present phases $\Q^n_K$ which is an additional discrete unknown of the system. This results from the fact that the molar fractions of mobile species and, in our case, the mineral concentrations corresponding to absent phases are removed from the set of unknowns. Regarding the set of equations in a given cell $K$, two questions arise if the phase of some primary species is absent: 
	\begin{itemize}
		\item[(i)] How to prevent the singularity of the total primary species molar conservation equations? 
		\item[(ii)] How to select the equilibrium and kinetic closure equations which must be kept in the system? 
	\end{itemize}
	To address these questions, given $\Q\subset \P$, we need to consider the restriction $E_{\Q}$ of the formula matrix $E$ to the set of present species $\alpha^{-1}(\Q)$ obtained by removing the columns corresponding to absent species $\alpha^{-1}(\Q^c)$: 
	$$
	E_{\Q} =
	\begin{pNiceArray}{c|c}[first-row,last-col]
		\E^\text{pr}\cap\alpha^{-1}(\Q) & \E^\text{sec}\cap\alpha^{-1}(\Q) &  \\
		\text{Id} & F_{\Q} & \E^\text{pr}\cap\alpha^{-1}(\Q) \\
		\hdashline
		0 & G_{\Q} & \E^\text{pr}\cap\alpha^{-1}(\Q^c)
	\end{pNiceArray} .
	$$
	The rank of $E_{\Q}$ is equal to
	$$
	\#\big(\E^\text{pr}\cap\alpha^{-1}(\Q)\big) + \rank G_{\Q}.
	$$
	To prevent the singularity of the total primary species molar conservation equations, we need to cope with the potential rank deficiency of $E_{\Q}$.
	This issue is handled by introducing, in each cell $K$, additional concentration unknowns
	$$
	\widetilde{n}_{e}, \quad e\in\widetilde{\E}_{\Q}, 
	$$
	where $\widetilde{\E}_{\Q}$ corresponds to a basis of a complementary subspace completing $\image(G_{\Q})$ in $\R^{\E^\text{pr}\cap \alpha^{-1}(\Q^c)}$.
	The number of such additional unknowns is
	$$
	\#\widetilde{\E}_{\Q} = \#\big(\E^\text{pr}\cap\alpha^{-1}(\Q^c)\big) - \rank G_{\Q}.
	$$
	It results that $\widetilde{\E}_{\Q}$ is such that the modified matrix 
	$$
	\widetilde E_{\Q} =
	\begin{pNiceArray}{c;c|c}[first-row,last-col]
		\E^\text{pr}\cap\alpha^{-1}(\Q) & \widetilde{\E}_{\Q}  & \E^\text{sec}\cap\alpha^{-1}(\Q) &  \\
		\text{Id} & 0 & F_{\Q} & \E^\text{pr}\cap\alpha^{-1}(\Q) \\
		\hdashline
		0 & \text{Id} & G^{(1)}_{\Q} & \widetilde{\E}_{\Q} \\
		\hdashline
		0 & 0 & G_{\Q}^{(2)} & \E^\text{pr}\cap\alpha^{-1}(\Q^c) \setminus \widetilde{\E}_{\Q}
	\end{pNiceArray} .
	$$
	is full rank, where we have applied a line reordering and denoted by $G^{(1)}_{\Q}$ and $G^{(2)}_{\Q}$ the submatrices of $G_{\Q}$ corresponding to the lines $\widetilde{\E}_{\Q}$ and $\E^\text{pr}\cap\alpha^{-1}(\Q^c) \setminus \widetilde{\E}_{\Q}$ respectively.
	This leads to the definition of the following set of cell unknowns depending on the set of present phases $\Q_K$ in cell $K$: 	
	\begin{equation}\label{dof_coats_cellK_QK}
		\begin{aligned}
			& P_K,\\
			& S^\beta_K,\quad \beta\in\P^\text{f},\\
			& \phi^\text{s}_{K},\\
			& C_{e,K},  \quad e\in \alpha^{-1}(\Q^\text{f}_K), \\
			& n_{m,K}, \quad m\in \alpha^{-1}(\Q^\text{s}_K),  \\
			& \widetilde{n}_{e,K}, \quad e \in\widetilde{\E}_{\Q_K}.  
		\end{aligned}
	\end{equation}
	
	Let us introduce the notation $\chi_K(e) := \1_{\Q_K}\big(\alpha(e)\big)$ which is one if the phase of $e$ is present in the cell $K$ and zero otherwise. Similarly, we denote by $\chi_K^{n-1}(e) := \1_{\Q_K^{n-1}}\big(\alpha(e)\big)$ the phase indicator at the previous time step.
	Then, the molar conservation of total primary species \eqref{eq:esp-base}, taking into account the present phases and the additional unknowns $(\widetilde{n}_{e,K})_{e \in\widetilde{\E}_{\Q_K}}$, is discretized as follows for each cell $K \in \cells_h$:
	\begin{multline}
		\1_{\widetilde{\E}_{\Q_K} }(e) |K| ~\widetilde{n}_{e,K} + \sum_{{e'}\in \E} E_{e,e'} |K| ~\Big(\chi_K(e')n_{e',K} - \chi_K^{n-1}(e') n^{n-1}_{e',K}\Big) \\
		+ \Delta t^n \sum_{{e'}\in \E^\text{f}} \ \sum_{\sigma=K|L \in \faces_K\cap \faces_h} \chi_{K_\sigma^{\alpha(e')}}(e') E_{e,e'} \Big(\zeta^{\alpha(e')} C_{e'} \eta^{\alpha(e')}\Big)_{K_\sigma^{\alpha(e')}} V_{KL}^{\alpha(e')} = 0, \quad e\in \E^\text{pr}. 
		\label{eq:reactive-coats}
	\end{multline}      
	The $M$ conservation equations \eqref{eq:reactive-coats} are coupled to the following local closure equations in each cell $K$:
	\begin{equation}\label{eq_fermeture_K}
		\begin{aligned}
			\sum_{\beta \in \Q^\text{f}_K} S^\beta_K  & = 1, & \\
			S^\beta_K  & = 0,  \quad \beta \in \P^\text{f} \setminus \Q^\text{f}_K,  \\
			\sum_{e\in \alpha^{-1}(\{\beta\})} C_{e}  &= 1, \quad \beta \in \Q^\text{f}_K,\\
			\phi^s_K - \sum_{m\in \Q^\text{s}_K} {n_{m,K}^s \over \zeta_{m,K}} & = 0. &   
		\end{aligned}  
	\end{equation}
	To close the system, it remains to specify the equations associated to chemical reactions. It is clear that, among the secondary species columns $e\in \E^\text{sec}\cap\alpha^{-1}(\Q_K)$ of the restricted formula matrix $E_{\Q_K}$, only the chemical reactions which do not involve absent primary species are kept in the system. These chemical reactions are obtained by linear combinations defined by the kernel of $G_{\Q_K}$. To be more specific, denoting by $v^{(i)}\in \R^{\E^\text{sec}\cap\alpha^{-1}(\Q_K)}$, for $i=1,\cdots, \dim\ker G_{\Q_K}$, a basis of $\ker G_{\Q_K}$, they are defined by 
	\begin{equation}\label{eq_kept_reaction}
		\sum_{e^\text{sec}\in \E^\text{sec}\cap\alpha^{-1}(\Q_K)} v^{(i)}_{e^\text{sec}} ~e^\text{sec} \leftrightarrows  \sum_{ e^\text{sec} \in \E^\text{sec}\cap\alpha^{-1}(\Q_K)} v^{(i)}_{e^\text{sec}} \sum_{e^\text{pr}  \in \E^\text{pr}\cap\alpha^{-1}(\Q_K)}   (F_{\Q_K})_{ e^\text{pr} ,e^\text{sec} } ~e^\text{pr}. 
	\end{equation}
	Hence, if \eqref{eq_kept_reaction} is a linear combination of equilibrium reactions, then the system is closed by writing the equilibrium equation associated to each of the chemical reactions \eqref{eq_kept_reaction}.
	
	When reaction \eqref{eq_kept_reaction} is a linear combination of one or more kinetic reactions, the definition of an associated kinetic law is not straightforward and requires further investigation, following the work of Hamlat \textit{et al} \cite{hamlat_mathematical_2020, hamlat_analysis_2021, erhel_projected_2024}. A common choice to circumvent this difficulty adopted in the numerical test cases presented in Section~\ref{sec:num} is to select primary species exclusively from the aqueous phase. In this setting, $\ker G_{\Q} = \{0\}$ for all $Q$, so no reaction of the form \eqref{eq_kept_reaction} occurs.
	
	By counting the number of equations \eqref{eq:reactive-coats}-\eqref{eq_fermeture_K}-\eqref{eq_kept_reaction} and the number of unknowns \eqref{dof_coats_cellK_QK} in the cell $K$, we obtain 
	$$
	\underbrace{\Big(1 + \#\P^\text{f} + 1 + \#\E - \#\alpha^{-1}(\Q_K^c) + \#\widetilde{\E}_{\Q_K}\Big)}_{\text{unknowns}} = \underbrace{\Big(1 + \#\P^\text{f} + 1 + \#\E^\text{pr} + \dim\ker G_{\Q_K} \Big)}_{\text{equations}},
	$$
	using, from the rank theorem, that
	$$
	\dim\ker G_{\Q_K} = \#\big(\E^\text{sec}\cap\alpha^{-1}(\Q_K)\big) - \rank G_{\Q_K}. 
	$$
	Finally, the set of present phases $\Q_K$ is updated at each Newton iteration for every cell $K$ in $\cells_h$ using phase appearance/disappearance criteria. As usual for the Coats' formulation,  fluid phase disappearance of a present phase is based on the negative sign of the phase saturation, while for a mineral phase its disapperance is based  on the negative sign of its concentration. Regarding phase appearance, it is typically based on negative flash calculations, stability analysis, or simplified approaches as discussed in the following numerical section. The additional unknowns $\widetilde{n}_{e,K}$, $e \in\widetilde{\E}_{\Q_K}$ not only avoid the singularity of the system, but they also track in \eqref{eq:reactive-coats} the possible appearance of absent primary species from neighboring cells $L$. Indeed, the sets $\widetilde{\E}_{\Q_K}$ are updated together with $\Q_K$ at each Newton iterate with  unknowns $\widetilde{n}_{e,K}$, $e\in \widetilde{\E}_{\Q_K}$ initialized to zero since they are associated to an absent phase.
        
	\section{Numerical results} \label{sec:num}
	
	The first objective of this section is to compare our numerical model to the results presented in the recent $\cod$ benchmark \cite{de_hoop_benchmark_2024, ahusborde_benchmark_2024}, both in terms of numerical solutions and computational efficiency. Then, the second objective is to assess the ability of our model formulation to simulate efficiently a $\hd$ storage test case in a 2D axisymmetric geometrical configuration.  	
	
	In both cases, the chemical system comprises two fluid phases, a liquid aqueous phase and a gas phase, denoted by $\P^\text{f} = \{\l,g\}$, and a certain number of mineral phases  $\P^\text{s} = \{m_1,\hdots, m_{N_\text{s}}\}$. The molar and mass densities of both fluid phases $\beta$ are related by 
	\begin{equation}
		\rho^\beta = \left(\sum_{e \in \alpha^{-1}(\{\beta\})} C_e M_e \right)\zeta^\beta,
		\label{eq:zeta-rho}
	\end{equation}
	with $M_e$ the molar mass of the species $e$.
	As will be specified in both test cases below, the liquid phase will be either slightly compressible or incompressible and the gas phase will either have a constant compressibility or behave as an ideal gas. The fluid phase relative permeabilities are modeled by the following Brooks-Corey law \cite{de_hoop_benchmark_2024}: 
	\begin{equation}
		k_\text{r}^\l = (S^\l)^2 \quad\text{and}\quad k_\text{r}^g = (S^g)^2.
		\label{eq:brooks-corey} 
	\end{equation}
	The absolute permeability tensor follows the Kozeny-Carman law \cite{hommel_porositypermeability_2018,de_hoop_benchmark_2024}: 
	\begin{equation}
		\mathbb{K} = \mathbb{K}_0 \left(\frac{\phi}{\phi_0}\right)^3,
		\label{eq:kozeny-carman}
	\end{equation}
	where $\phi_0$ and $\mathbb{K}_0$ are initial porosity and permeability respectively.
	
	Species $e$ involved in liquid gas phase transitions are denoted as a pair defined by their component $c$ in $\C$, the set of all components, and their phase $\beta$ in $\P^\text{f}$ i.e. we set $e = (c,\beta)$. When both liquid and gas phases are present, the phase equilibrium  is written in terms of partition coefficient $K_c$ as follows for all components  $c \in\C$. 
	\begin{equation}
		C_{(c,g)} = K_{c} \, C_{(c,l)}.
		\label{eq:phase-equi}
	\end{equation}
	Each other equilibrium reactions is modeled by the \textit{mass action law} associated with its secondary species $e^\text{sec}$ linked to its reaction $R \in \Rea$, namely (see \cite[Section~4.1.2]{de_hoop_benchmark_2024})
	\begin{equation}
		\Upsilon := \frac{1}{a_{e^\text{sec}}}\prod_{e^{\text{pr}} \in \E^{\text{pr}}} \big( a_{e^{\text{pr}}} \big)^{E_{e^{\text{pr}},R}} = K_R,
		\label{eq:mass-action-law}
	\end{equation}
	where $a_e$ is the activity of $e$ and $K_R$ [$-$] is the equilibrium constant associated to the reaction $R$.
	
	If kinetic reactions are considered instead, the conservation equation associated to the secondary species is retained in its original non-Morel form and its reaction rate is modeled through appropriate kinetic rate laws detailed for each test case below.
	
	At each time step $n$, the nonlinear system is solved using a Newton method. If the Newton iteration converges, the computation proceeds to the next time step and the time-step size $\Delta t^{n+1}$ is increased, up to a prescribed maximum value. More precisely, for a fixed $\Delta t_\text{max}$ and $\Delta t^{1}$ given
	$$ \Delta t^{n+1} = \min\{1.1\, \Delta t^n, \Delta t_\text{max}\}, $$ 
	if the Newton iteration converges for the time step $n$.
	If the Newton method fails to converge, the time-step size is halved and the current time step is recomputed. If the time-step size becomes too small, then the simulation fails.
	
	At each Newton iteration, the Jacobian matrix and residual vector are assembled, and a relaxed Newton step is computed using a weighted infinity norm based on pressure, saturations, fluid-species compositions, and mineral molar concentrations. The nonlinear iteration is considered converged when either the infinity norm of the relative residual drops below the prescribed tolerance ($10^{-7}$), or the infinity norm of the scaled Newton increment becomes smaller than $10^{-10}$.
	
	To solve the linear system arising from the Newton iterations, the CPR-AMG solver is employed. It combines the CPR (Constraint Pressure Residual) preconditioner with an Algebraic Multigrid (AMG) solver \cite{cao_parallel_2005}. The linear solver converges when the relative two-norm of the residual becomes sufficiently small, typically smaller than $10^{-7}$.
	
	The implementation of the algorithm is done in Fortran.   
	
	\subsection{Phase appearance and disappearance criteria} \label{sec:phase-dis}
	
	A phase is removed when it becomes nonphysical after the Newton update. More precisely, a fluid phase disappears if its saturation becomes negative, while a mineral phase $m$ is removed if its concentration $n_m$ becomes negative. We note that in \eqref{eq:reactive-coats}, the flux $V_{KL}^\beta$ does not need to be computed if phase $\beta$ is absent from both cells $K$ and $L$, since in that case the mobility terms at face $\sigma=K|L$ vanish. The density in the gravity term of the Darcy flux is therefore always well defined by a saturation weighted average between both cell values.
	
	The appearance of new phases is then evaluated. When the liquid phase is present, an absent mineral phase is introduced if its solubility product exceeds the equilibrium constant. Similarly, a gas phase appears if $\widetilde{n}_e > 0$ for some gas species $e \in \widetilde{\E}_{Q_K}$, or if
	$$ \sum_{c\in \C} \widetilde{C}_{(c,g)} > 1, $$
	where $\widetilde{C}_{(c,g)}$ denotes the gaseous mole fractions at equilibrium with the liquid phase i.e. that satisfies \eqref{eq:phase-equi} which provides the gas phase appearance criteria: 	
	$$ \sum_{ c\in \C} K_c \, C_{(c,\l)} > 1. $$	
	When the liquid phase is absent, phase appearance is determined by performing a multiphase (negative) flash calculation \cite{leal_overview_2017} at fixed total composition. For every primary species $e \in \E^\text{pr}$, we define its total composition
	$$ z_e := \sum_{\beta \in \P^{\mathrm{f}}} \theta^\beta \sum_{{e'}\in \E^\text{f}} E_{e,e'} C_{e'} + \sum_{{e'}\in \E^\text{s}} E_{e,e'} n_{e'}. $$
	We include phase change equation \eqref{eq:phase-equi} as well as the mass action law for the chemical reactions \eqref{eq:mass-action-law} in the fluid phases. Here, we only consider equilibrium reactions. Finally, the system is closed using \eqref{eq:close-comp-fluid}. To reduce the nonlinearity of the system, we introduce the composition parametrization proposed in~\cite{jonval_reactive_2025}.
	
	The system of $(\#\E + \#\P^\text{f})$ equations is solved in terms of fluid compositions, of mineral concentrations and of $\theta^\beta$ [$\mathrm{mol \cdot m^{-3}}$] for $\beta$ a fluid phase (i.e. $(\#\E + \#\P^\text{f})$ unknowns) in each cell without liquid phase. Since the system is nonlinear, we use a Newton method whose initialization is chosen as, for all $e^\text{f}\in\E^\text{f}$ in phase $\beta$ and $e^\text{s}\in\E^\text{s}$
	$$ C_{e^\text{f}}^{0} = \frac{1}{\#\alpha^{-1}(\{\beta\})}, \quad n_{e^\text{s}}^{0} = n_{e^\text{s}} \quad\text{and}\quad \theta^{\beta,0} = \phi \zeta^\beta S^\beta, $$
	with $n_{e^\text{s}}$, $\phi$, $\zeta^\beta$ and $S^\beta$ from the current Newton iteration of the implicit scheme. Finally, to improve the convergence of the Newton method, we implement a backtracking algorithm using the Armijo criterion \cite{press_numerical_2007}. Then, the liquid phase is introduced if its computed mole fraction $\theta^\l$ is positive. 
	
	When a fluid phase appears, it is common to initialize its saturation to zero, which may lead to a singular Jacobian matrix. One possible remedy is to regularize the Jacobian. First, we detect the vanishing Jacobian columns associated to an unknown $\delta C_{e,K}$, where $e$ is a primary species. Then, the linearized conservation equation of the primary species $e$ in cell $K$ is replaced by the equation $\delta C_{e,K} = 0$. Another approach consists in replacing the mole fraction unknowns of this phase by the variable $S^\beta C_e$, $e$ in phase $\beta$. Lastly, a more straightforward alternative, adopted in the sequel, is to initialize the appearing phase saturation with a small non-zero value.
	
	\subsection{$\cod$ benchmark test case}
	
\begin{figure}[H]
		\centering
		\begin{tikzpicture}[x=0.015cm, y=0.015cm, >=stealth, font=\sffamily]
			
			\definecolor{dom1}{RGB}{175, 198, 233}      
			\definecolor{dom2}{RGB}{191, 191, 191}      
			\definecolor{gasgreen}{RGB}{0, 176, 80}     
			\definecolor{waterblue}{RGB}{0, 112, 192}   
			\definecolor{outred}{RGB}{255, 0, 0}        
			
			\draw[thick, fill=dom2] (0, 0) rectangle (600, 240);
			
			\draw[thick, fill=dom1] (120, 80) rectangle (480, 160);
			
			\node[anchor=south west] at (5, 5) {$\Omega_2: \{\phi_2, k_2\}$};
			\node[anchor=south west] at (125, 85) {$\Omega_1: \{\phi_1, k_1\}$};
			
			
			\draw[thick] (0, -30) -- (600, -30);
			\foreach \x in {0, 120, 480, 600} \draw[thick] (\x, -15) -- (\x, -45);
			\node[fill=white] at (60, -30) {120};
			\node[fill=white] at (300, -30) {360};
			\node[fill=white] at (540, -30) {120};
			
			\draw[thick] (-60, 0) -- (-60, 240);
			\foreach \y in {0, 120, 240} \draw[thick] (-45, \y) -- (-75, \y);
			\node[rotate=90, fill=white, inner sep=2pt] at (-60, 60) {120};
			\node[rotate=90, fill=white, inner sep=2pt] at (-60, 180) {120};
			
			\draw[thick] (660, 0) -- (660, 240);
			\foreach \y in {0, 80, 160, 240} \draw[thick] (645, \y) -- (675, \y);
			\node[rotate=90, fill=white, inner sep=2pt] at (660, 40) {80};
			\node[rotate=90, fill=white, inner sep=2pt] at (660, 120) {80};
			\node[rotate=90, fill=white, inner sep=2pt] at (660, 200) {80};

			\tikzset{streamarrow/.style={-{Triangle[length=2mm, width=2mm]}, line width=2pt}}
			
			\foreach \y in {15, 45, 75, 105} {
				\draw[streamarrow, gasgreen] (-45, \y) -- (-5, \y);
			}
			
			\foreach \y in {135, 165, 195, 225} {
				\draw[streamarrow, waterblue] (-45, \y) -- (-5, \y);
			}
			
			\tikzset{outputarrow/.style={-{Triangle[length=2mm, width=2mm]}, line width=2pt}}
			
			\foreach \y in {12, 36, 60, 84, 108, 132, 156, 180, 204, 228} {
				\pgfmathsetmacro{\xend}{605 + 35 - (\y / 240) * 20}
				\draw[outputarrow, outred] (605, \y) -- (\xend, \y);
			}		
			
			\begin{scope}[shift={(60, -100)}]
				\draw[-{Triangle[length=2mm, width=2mm]}, line width=2pt, gasgreen]  (0, 0) -- (40, 0);
				\node[anchor=west] at (50, 0) {Gas stream (pure $\codg$):};
				\node[anchor=west] at (390, 0) {Q$_{\text{gas}}$ [m$^3$/day]};
				
				\draw[-{Triangle[length=2mm, width=2mm]}, line width=2pt, waterblue] (0, -35) -- (40, -35);
				\node[anchor=west] at (50, -35) {Water stream (pure $\hdol$):};
				\node[anchor=west] at (390, -35) {Q$_{\text{wat}}$ [m$^3$/day]};
				
				\draw[-{Triangle[length=2mm, width=2mm]}, line width=2pt, outred]    (0, -70) -- (40, -70);
				\node[anchor=west] at (50, -70) {Outflow (hydrostatic pressure):};
				\node[anchor=west] at (390, -70) {p$_{\text{out}}$ [bar]};
			\end{scope}
			
		\end{tikzpicture}
		\caption{Configuration of the $\cod$ injection test case.}
		\label{fig:config-cod-ben}
	\end{figure}
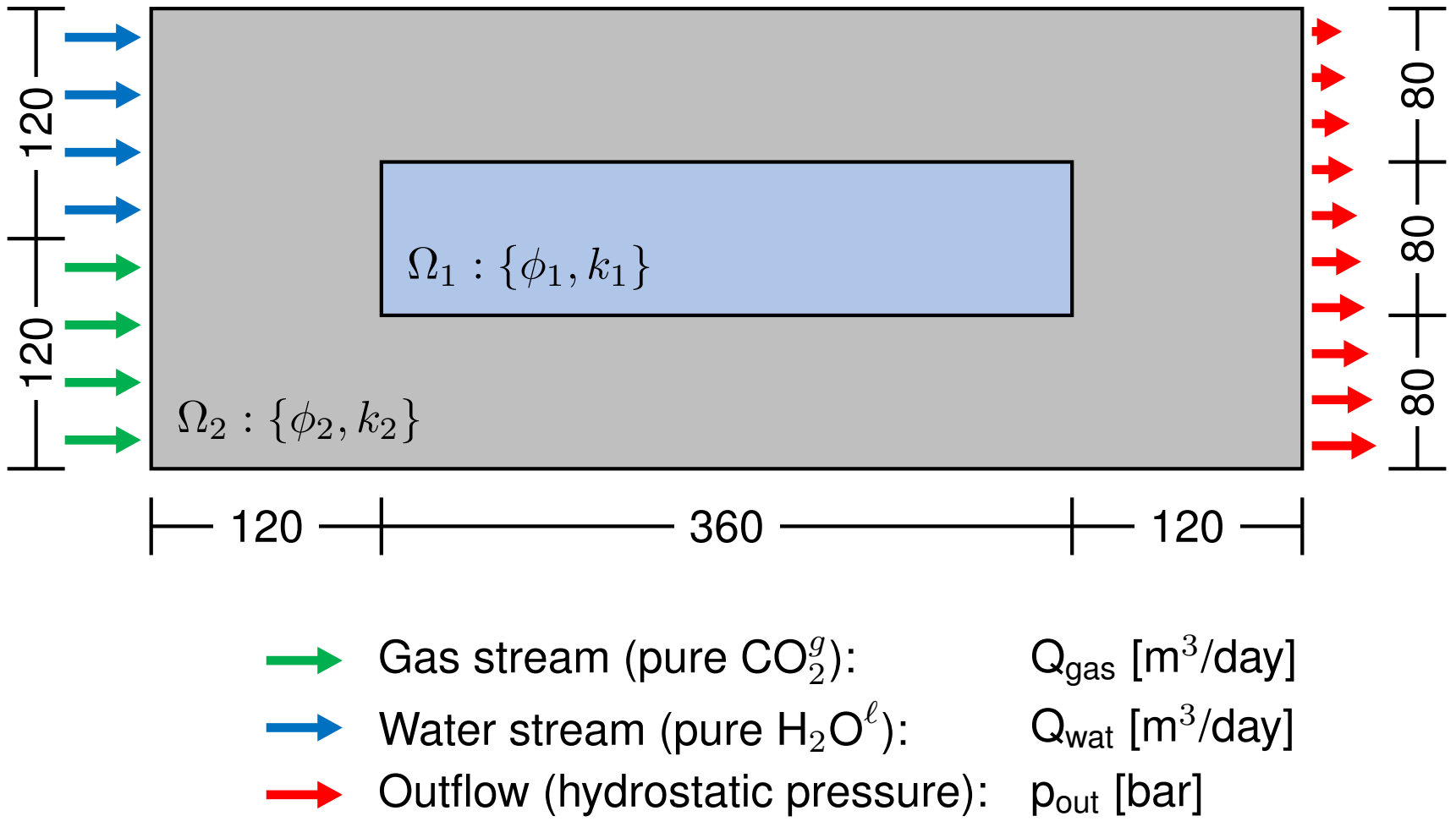

	We compare our results with the two-dimensional \textit{extended chemical model} with gravity\footnote{available at \url{https://github.com/eahusbor/Reactive-Multiphase-Benchmark}} described in \cite[Section~{5.2.2}]{de_hoop_benchmark_2024} and in \cite{ahusborde_benchmark_2024}.
	For the reader's convenience, we recall below the key features of the dataset and refer to  \cite{de_hoop_benchmark_2024} for its complete description. This test case, illustrated in Figure~\ref{fig:config-cod-ben}, considers a reservoir of dimensions $240\,\mathrm{m} \times 600\,\mathrm{m}$, composed of two different rock types denoted by $\Omega_1$ and $\Omega_2$. Liquid water is injected at a constant flux $Q_\text{wat}$ through the upper half of the left boundary, while gaseous $\cod$ is injected at a constant flux $Q_\text{gas}$ through the lower half. We consider Dirichlet boundary conditions that account for hydrostatic pressure on the right boundary. The reservoir is initially assumed to be in chemical equilibrium.
	The fluid phases are characterized by their mass density, for $\beta\in\{\l,g\}$
	\begin{equation}
		\rho^\beta(P) = \rho^\beta_0 \big(1 + c^\beta(P - P_0)\big),
		\label{eq:lin-compress-fluid}
	\end{equation}
	with $\rho^\beta_0$ the reference density at pressure $P_0$ set to $100$~$\mathrm{kg \cdot m^{-3}}$ for the gaseous phase and to $1000$~$\mathrm{kg \cdot m^{-3}}$ for the aqueous phase, and $c^\beta$ its compressibility [$\mathrm{Pa^{-1}}$] fixed at $10^{-9}$ for the gaseous phase and at $10^{-11}$ for the aqueous phase.
	Aqueous species activities are expressed in terms of molality as
	$$ a_e = M_{\hdo^\l} \frac{C_e}{C_{\hdo^\l}},\quad \forall e\in\alpha^{-1}(\{\l\}), $$
	where $M_{\hdo^\l}$ denotes the number of moles of $\hdo$ per kilogram of aqueous phase (typically taken as $55.508$~$\mathrm{mol \cdot kg^{-1}}$).	When the liquid phase disappears, since $C_{\hdo^\l}$ tends to zero, which leads to a singular behavior in the activity definition. To overcome this issue, we regularize the liquid-phase activities as
	$$ a_e = M_{\hdo^\l} \frac{C_e}{C_{\hdo^\l} + \varepsilon},\quad \forall e\in\alpha^{-1}(\{\l\}), $$
	where $\varepsilon > 0$ is a small regularization parameter, set to $10^{-15}$. The reaction rate $T_{\cacot}$ for the $\cacot$ kinetic reaction is given in \cite{de_hoop_benchmark_2024} and reads
	$$ T_{\cacot} = K_{\cacot}\,(1 - \phi_0) \left(1 - \frac{\Upsilon}{K_{\cacot}^\text{sp}}\right)n_{\cacot}. $$
	where $K_{\cacot}$ [$\mathrm{mol \cdot m^{-2} \cdot s^{-1}}$] is the kinetic rate constant, $\Upsilon$ is the ion activity product \eqref{eq:mass-action-law}, and $K_{\cacot}^\text{sp}$ is the equilibrium constant. 	
	It is assumed that the liquid phase follow a linear compressibility model.
	The constants and the initial conditions are given in \cite{de_hoop_benchmark_2024}. 
	
	Using Morel's formalism described in Section~\ref{sec:model}, we choose for this test case the following set of primary species:
	$$ \E^\text{pr} = \{ \hdol, \codl, \hp, \cadp \}, $$
	and of secondary species:
	$$ \E^\text{sec} = \{ \hcotm, \ohm, \cotdm, \hdog, \codg, \cacot \}. $$
	The six reactions in $\Rea$ reads
	$$ \begin{aligned}
		\hcotm &  \leftrightarrows \hdol + \codl - \hp, & \text{(equilibrium aq-aq)}\\
		\ohm & \leftrightarrows \hdol - \hp, & \text{(equilibrium aq-aq)} \\
		\cotdm & \leftrightarrows \hdol + \codl - 2\hp, & \text{(equilibrium aq-aq)}\\
		\hdog & \leftrightarrows \hdol, & \text{(phase change)} \\
		\codg & \leftrightarrows \codl, & \text{(phase change)} \\
		\cacot & \leftrightarrows \cadp + \hdol + \codl -2 \hp, & \text{(kinetic min-aq)}  
	\end{aligned} $$
	and the set of phases $\P$ is defined by 
	$$ \P = \{ \l, g, m_{\cacot} \}. $$
	This leads to the following formula matrix,
	$$ E := \begin{pNiceArray}{cccc|cccccc}[first-row,last-col]
		\hdol & \codl & \hp & \cadp & \hcotm & \ohm & \cotdm & \hdog & \codg & \cacot \\
		1 & 0 & 0 & 0 & 1 & 1  & 1  & 1 & 0 & 1 & \hdol \\
		0 & 1 & 0 & 0 & 1 & 0  & 1  & 0 & 1 & 1 & \codl \\
		0 & 0 & 1 & 0 &-1 & -1 & -2 & 0 & 0 & -2 & \hp \\
		0 & 0 & 0 & 1 & 0 & 0  & 0  & 0 & 0 & 1 & \cadp
	\end{pNiceArray} .$$
	Finally, we need to deal with the singularity of the system of equations when the aqueous phase disappears for $\Q_K^\text{f} = \{g\}$.
	Indeed, in that case no primary species remain as they all belong to the liquid phase. Two subcases need to be considered, first if $\Q_K^\text{s}=\emptyset$, then
	$$ E_{\Q_K} = \begin{pNiceArray}{cc}[first-row,last-col]
		\hdog & \codg \\
		1 & 0 & \hdol \\
		0 & 1 & \codl \\
		0 & 0 & \hp \\
		0 & 0 & \cadp
	\end{pNiceArray}, $$
	meaning that we need to add the two additional unknowns $\widetilde{n}_\hp$ and $\widetilde{n}_\cadp$. Second, if $ \Q_K^\text{s}=\{ m_{\cacot} \}$, then
	$$ E_{\Q_K} = \begin{pNiceArray}{ccc}[first-row,last-col]
		\hdog & \codg & \cacot \\
		1 & 0 & 1 & \hdol \\
		0 & 1 & 1 & \codl \\
		0 & 0 & -2 & \hp \\
		0 & 0 & 1 & \cadp
	\end{pNiceArray}, $$
	which requires adding one additional unknown  $\widetilde{n}_\hp$ (or $\widetilde{n}_\cadp$).
	To sum up, we choose
	$$ \widetilde{\E}_{\Q_K} = \begin{cases}
		\{ \hp, \cadp \} & \text{if}\ \Q_K^\text{s}=\emptyset \text{ and } \Q_K^\text{f} = \{g\},\\
		\{ \hp \} & \text{if}\ \Q_K^\text{s}=\{ m_{\cacot} \} \text{ and } \Q_K^\text{f} = \{g\}.
	\end{cases} $$
	Since all primary species are absent if the aqueous phase is absent, it is clear that no reaction remains for $\Q_K^\text{f} = \{g\}$.

	Following \cite{de_hoop_benchmark_2024}, the simulation is performed on a uniform Cartesian mesh of size of $120\times 48$.
	Table~\ref{tab:benchmark-cod} reports the total number of time steps performed during the simulation, the average number of nonlinear iterations per time step, the average number of linear iterations per Newton iteration, the number of flash computations and the CPU time for DuMu$^\text{X}$, CooresFlow and our simulations (Coats-TR). The number of failed attempts for each quantity is indicated in parentheses.
	Our performance is comparable to that of CooresFlow, while remaining faster, mainly due to the use of a more recent CPU architecture. It should also be noted that flash computations are performed locally in each grid cell and could be easily parallelized across multiple processes, although this is not implemented in the current version.
	
\begin{figure}[H]
		\subfloat{
			\includegraphics[width=\2dplotw]{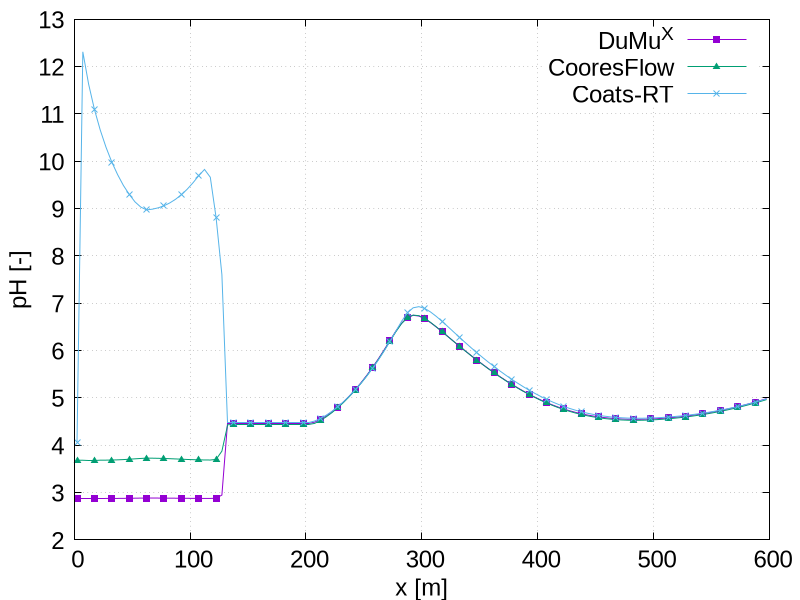}
		}
		\subfloat{
			\includegraphics[width=\2dplotw]{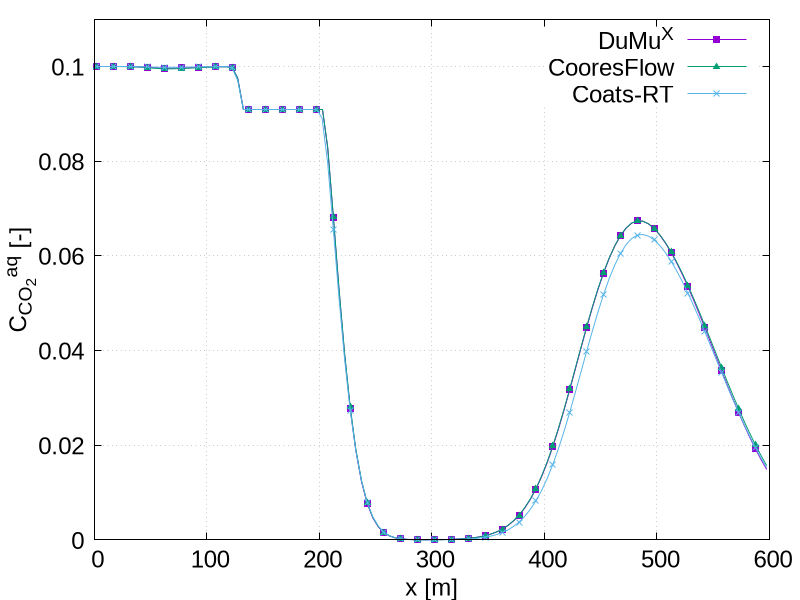}
		}
		\caption{$\cod$ benchmark test case: Comparison of the cross-section at $y = 50$ m of gas saturation, pressure, pH and molar fraction of gaseous $\cod$ obtained using the DuMu$^\text{X}$ and CooresFlow simulation codes and our method.}
		\label{fig:ben-cod-plot}
	\end{figure}

	Figure~\ref{fig:ben-cod-plot} shows a comparison of the cross-section at $y = 50$ m obtained using the DuMu$^\text{X}$ and CooresFlow simulation codes and our method at time $t=1000$ days.
	A good agreement with the DuMu$^\text{X}$ and CooresFlow benchmarks is observed overall, with only slight deviations at the peak values of $C_{\codl}$ and pH.
	Regarding the pH, in regions where the liquid phase is absent, namely for $0 \leq x \leq 120$, the computed values are not physically meaningful.
	
\begin{figure}[H]
		\centering
		\subfloat[gas saturation $S^g$ at $t=500$ days]{
			\includegraphics[width=\2dimgw]{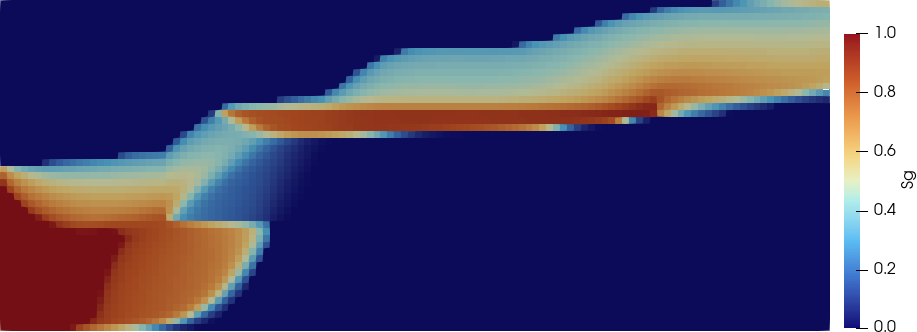}
		}
		\subfloat[gas saturation $S^g$ at $t=1000$ days]{
			\includegraphics[width=\2dimgw]{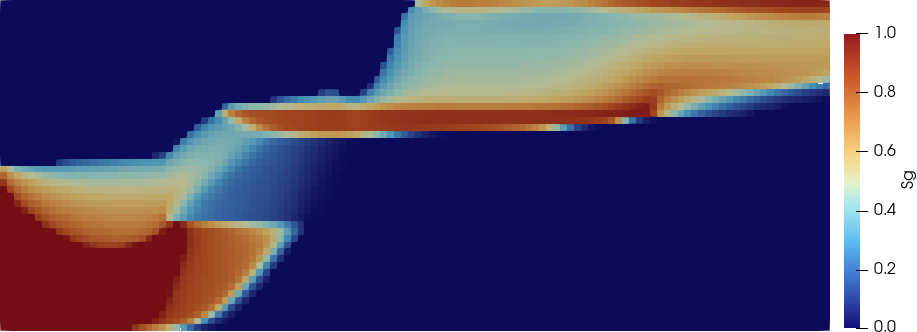}
		}
		\quad
		\subfloat[molar fraction $C_{\codg}$ at $t=500$ days]{
			\includegraphics[width=\2dimgw]{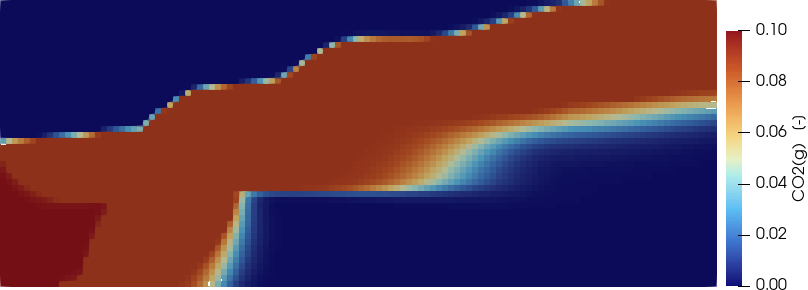}
		}
		\subfloat[molar fraction $C_{\codg}$ at $t=1000$ days]{
			\includegraphics[width=\2dimgw]{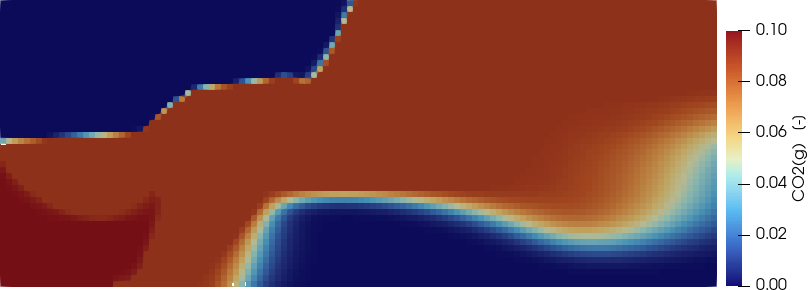}
		}
		\caption{$\cod$ benchmark test case: Gas saturation and the molar fraction of gaseous $\cod$ at times $t = 500$ days and $t = 1000$ days using our method.}
		\label{fig:ben-cod-2d}
	\end{figure}

	Figure~\ref{fig:ben-cod-2d} presents the gas saturation and the molar fraction of gaseous $\cod$ at times $t = 500$ days and $t = 1000$ days. The results are consistent with those obtained using DuMu$^\text{X}$ and CooresFlow.
	
	\begin{table}[H]
		\centering
		\begin{tabular}{|c|c|c|c|c|c|}
			\hline
			\rowcolor{gray!40} Code team & TS (failed) & NI & LI / NI & F (failed) & CPU \\
			\hline
			DuMu$^\text{X}$ & 17571 (2937) & 3.45 & 8.37 & - & 66124  \\
			CooresFlow & 2008 (0) & 4.1 & not available & - & 10224 \\
			Coats-RT & 2032 (8) & 5.2 & 11.76 & 206939 (0) & 6361 \\
			\hline
		\end{tabular}\vspace{0.2cm}		
		\caption{Numerical performance of the codes for $\cod$ benchmark. TS: number of time steps, NI: Average number of nonlinear iterations per time step, LI / NI: average number of linear iterations per Newton iteration, F: number of flash computations, CPU: elapsed time in seconds. Used architecture for DuMu$^\text{X}$: Intel i7-8565U @1.8 GHz, for CooresFlow: Intel Xeon E5-1620 v3 @3.5 GHz and for Coats-RT: Intel Core Ultra 7 255H @2.0 GHz}
		\label{tab:benchmark-cod}
	\end{table}
	
	\def\2dplotw{7.5cm}
	\begin{figure}[H]
		\centering
		\subfloat{
			\includegraphics[width=\2dplotw]{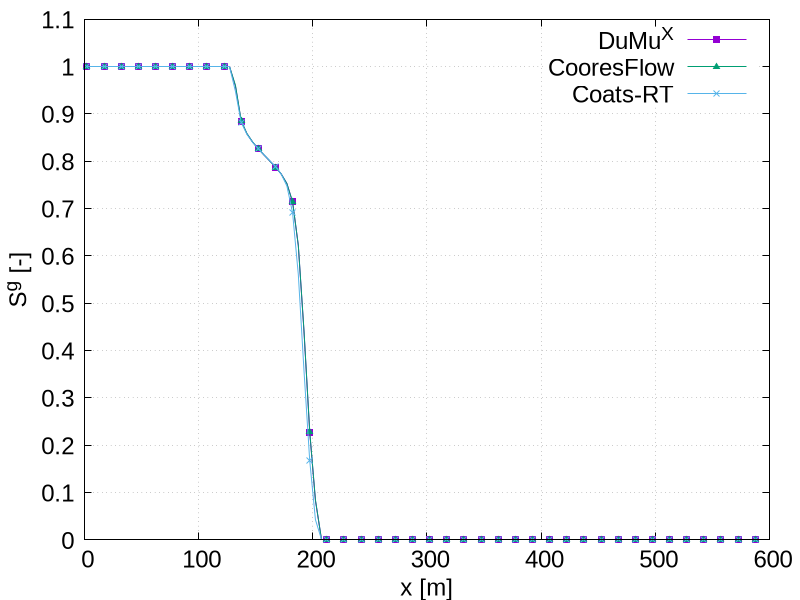}
		}
		\subfloat{
			\includegraphics[width=\2dplotw]{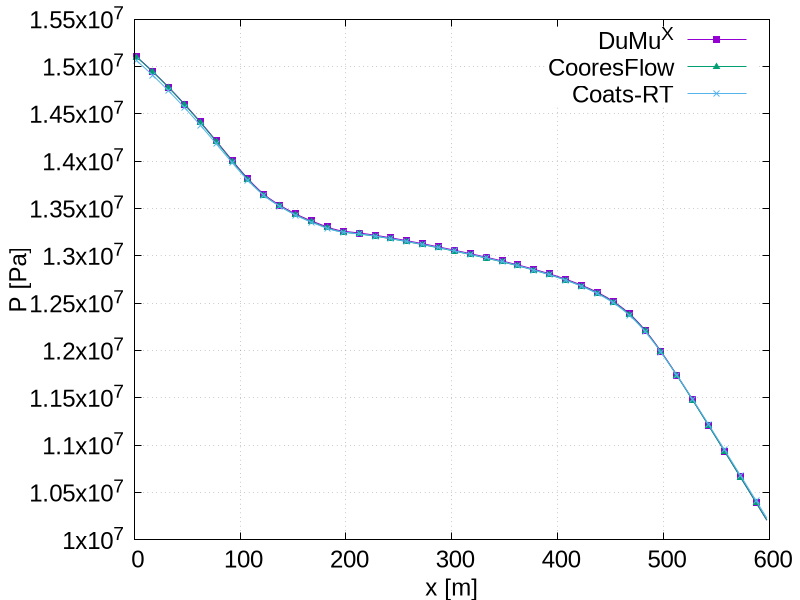}
		}
	\end{figure}
	
	\def\2dimgw{8cm}
	
	\subsection{$\hd$ storage test case}
	
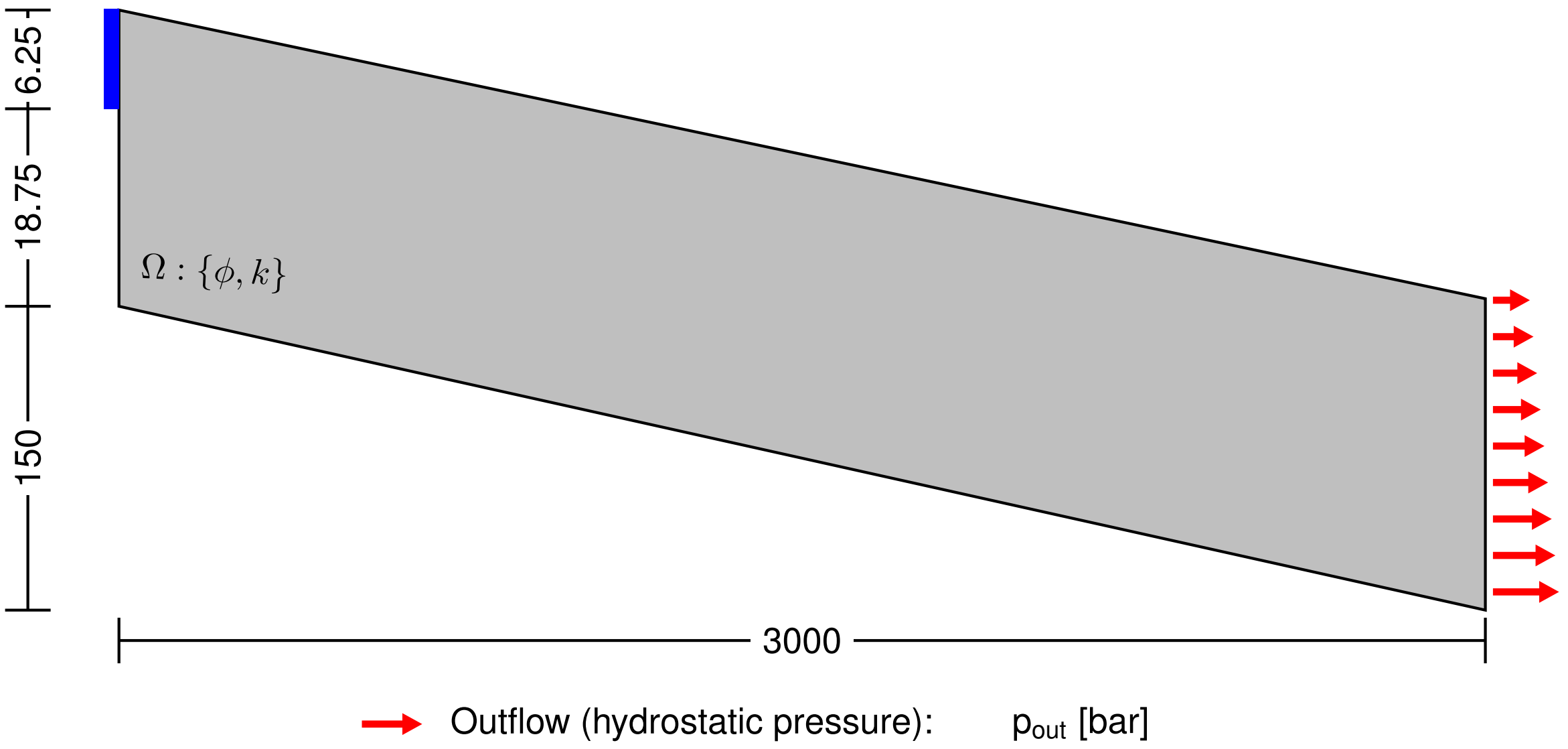
\begin{figure}[H]
		\centering
		\begin{tikzpicture}[scale=1, x=0.015cm, y=0.015cm, >=stealth, font=\sffamily]
			
			\definecolor{dom1}{RGB}{175, 198, 233}      
			\definecolor{dom2}{RGB}{191, 191, 191}      
			\definecolor{gasgreen}{RGB}{0, 176, 80}     
			\definecolor{waterblue}{RGB}{0, 112, 192}   
			\definecolor{outred}{RGB}{255, 0, 0}        
			
			\draw[thick, fill=dom2] (0, 0) -- (900, -200) -- (900, 05) -- (0, 195) -- cycle;
			
			\node[anchor=south west, rotate=-2.86] at (5, 5) {$\Omega: \{\phi, k\}$};
			
			
			\draw[thick] (0, -220) -- (900, -220);
			\draw[thick] (0, -205) -- (0, -235);
			\draw[thick] (900, -205) -- (900, -235);
			\node[fill=white] at (450, -220) {3000};
			
			\draw[thick] (-60, 0) -- (-60, 195);
			\foreach \y in {0, 130} \draw[thick] (-45, \y) -- (-75, \y);
			\node[rotate=90, fill=white, inner sep=2pt] at (-60, 65) {18.75};            
			
			\draw[thick] (-60, -200) -- (-60, 0);     
			\draw[thick] (-45, -200) -- (-75, -200);
			\node[rotate=90, fill=white, inner sep=2pt] at (-60, -100) {150};

			\draw[thick] (-60, 140) -- (-60, 195);
			\draw[thick] (-45, 195) -- (-75, 195);
			\node[rotate=90, fill=white, inner sep=2pt] at (-60, 162) {6.25};
			
			\tikzset{streamarrow/.style={-{Triangle[length=2mm, width=2mm]}, line width=2pt}}
			
			\tikzset{outputarrow/.style={-{Triangle[length=2mm, width=2mm]}, line width=2pt}}

			\draw[draw=none, fill=blue] (0, 130) rectangle (-10, 196);
			
			\foreach \h in {12, 36, 60, 84, 108, 132, 156, 180, 204} {
				\pgfmathsetmacro{\y}{\h-90 - 30-80}
				\pgfmathsetmacro{\xend}{905 + 45 - (\h / 200) * 20}
				\draw[outputarrow, outred] (905, \y) -- (\xend, \y);
			}       
			
			\begin{scope}[shift={(160, -260)}]
				
				\draw[-{Triangle[length=2mm, width=2mm]}, line width=2pt, outred]    (0, -15) -- (40, -15);
				\node[anchor=west] at (50, -15) {Outflow (hydrostatic pressure):};
				\node[anchor=west] at (420, -15) {p$_{\text{out}}$ [bar]};
			\end{scope}
			
		\end{tikzpicture}
		\caption{Radial cross-section of the axisymmetric reservoir for the $\hd$ storage test case. The well boundary located at radius $r_\mathrm{w}$ is represented in blue and the right boundary is located at radius $r=r_\mathrm{w}+3000$ m.}
		\label{fig:config-hd-ben}
	\end{figure}

	Figure~\ref{fig:config-hd-ben} depicts the setup for the $\hd$ storage test case in an anticline-shaped reservoir represented by an axisymmetric model based on an inclined radial cross-section of size $3000 \text{ m} \times 25 \text{ m}$. Starting from an initial chemical equilibrium state (given in Table~\ref{tab:hd-init} in Appendix~\ref{sec:values-h2}) at initial pressure $P_0$ in the reservoir, gaseous $\text{H}_2$ is injected through well on the top-left boundary (in blue on Figure~\ref{fig:config-hd-ben}) at a prescribed total molar rate, while the right boundary is maintained at a fixed hydrostatic outflow pressure accounting for gravity. The liquid phase is assumed to be incompressible with fixed mass density, while the gas phase behaves like an ideal gas with molar density defined, at fixed temperature $T_0=25$°C, by
	\begin{equation}
		\zeta^g(P) = \frac{P}{R\,T_0},
		\label{eq:ideal-gaz}
	\end{equation}
	with $R=8.314462618$ $[\mathrm{J}\cdot\mathrm{mol}^{-1}\cdot\mathrm{K}^{-1}]$ the ideal gas constant. The viscosity of the gaseous phase is taken as the viscosity of $\hd$, since the phase is composed almost entirely of hydrogen. The simulation then proceeds through successive cycles of $\hdg$ injection and production via a vertical well of radius $r_\mathrm{w}$ located at the center of the domain. The (simplified) chemical system includes a kinetic $\text{H}_2$-based sulfate reduction (representing reservoir microbial activity), aqueous acid-base reactions, and calcite dissolution/precipitation.
	This setup is designed to assess the code's robustness and efficiency when handling a significantly larger number of chemical species and reactions than for the $\text{CO}_2$ benchmark test case. The parameter values of the test case are given in Table~\ref{tab:param-h2} in Appendix~\ref{sec:values-h2}.
	
	\subsubsection{Chemical System}
	
	We choose the following set of primary species:
	$$ \E^\text{pr} := \{ \hdsl, \hdl, \codl, \hdo, \hp, \cadp, \nap \}, $$
	and of secondary species:
	$$ \E^\text{sec} = \{ \soqdm, \hsm, \ohm, \hcotm, \cotdm, \hdsg, \hdg, \codg, \cacot \}. $$
	The reactions in $\Rea$ reads
	$$ \begin{aligned}
		\soqdm & \leftrightarrows \hdsl - 4\hdl + 4\hdo - 2\hp, & \text{(kinetic aq-aq)} \\
		\hsm & \leftrightarrows \hdsl - \hp, & \text{(equilibrium aq-aq)} \\
		\ohm & \leftrightarrows \hdo - \hp, & \text{(equilibrium aq-aq)} \\
		\hcotm & \leftrightarrows \codl - \hp + \hdo, & \text{(equilibrium aq-aq)} \\
		\cotdm & \leftrightarrows \codl - 2\hp + \hdo, & \text{(equilibrium aq-aq)} \\
		\hdsg & \leftrightarrows \hdsl, & \text{(phase change)} \\
		\hdg & \leftrightarrows \hdl, & \text{(phase change)} \\
		\codg & \leftrightarrows \codl, & \text{(phase change)} \\
		\cacot & \leftrightarrows \cadp + \codl - 2\hp + \hdo, & \text{(equilibrium min-aq)}
	\end{aligned} $$ 
	and the the set of phases $\P$ is defined by 
	$$ \P = \{\l, g, m_{\cacot}\}. $$
	The associated formula matrix is
	$$ E := \left(\begin{array}{ccccccc|ccccccccc}
		1 & 0 & 0 & 0 & 0 & 0 & 0 &  1 &  1 &  0 &  0 &  0 & 1 & 0 & 0 &  0 \\
		0 & 1 & 0 & 0 & 0 & 0 & 0 & -4 &  0 &  0 &  0 &  0 & 0 & 1 & 0 &  0 \\
		0 & 0 & 1 & 0 & 0 & 0 & 0 &  0 &  0 &  0 &  1 &  1 & 0 & 0 & 1 &  1 \\
		0 & 0 & 0 & 1 & 0 & 0 & 0 &  4 &  0 &  1 &  1 &  1 & 0 & 0 & 0 &  1 \\
		0 & 0 & 0 & 0 & 1 & 0 & 0 & -2 & -1 & -1 & -1 & -2 & 0 & 0 & 0 & -2 \\
		0 & 0 & 0 & 0 & 0 & 1 & 0 & 0 &  0 &  0 &  0 &  0 & 0 & 0 & 0 &  1 \\
		0 & 0 & 0 & 0 & 0 & 0 & 1 & 0 &  0 &  0 &  0 &  0 & 0 & 0 & 0 &  0
	\end{array}\right). $$
	
	Note that $\nap$ is added to ensure electro-neutrality of the solution.
	
	To handle potential liquid phase disappearance, we set:  
	$$ \widetilde{\E}_{\Q_K} = \begin{cases}
		\{ \hdo, \hp, \cadp, \nap \} & \text{if}\ \Q_K^\text{s}=\emptyset \text{ and } \Q_K^\text{f} = \{g\},\\
		\{ \hdo, \hp, \nap \} & \text{if}\ \Q_K^\text{s}=\{ m_{\cacot} \} \text{ and } \Q_K^\text{f} = \{g\}.
	\end{cases} $$
	As in the previous test case, since all primary species are absent if the aqueous phase is absent, no reaction remains for $\Q_K^\text{f} = \{g\}$.
	
	Liquid species activities are defined as molar concentrations: for $e$ in $\alpha^{-1}(\{\l\})$
	$$ a_e = 10^{-3}\zeta^\l C_e. $$
	The aqueous homogeneous kinetic reaction follows a Monod-type law \cite{monod_technique_1978}
	$$ T_{\soqdm} = -v_\text{max}\, \frac{10^{-3}\zeta^\l C_{\soqdm}}{10^{-3}\zeta^\l C_{\soqdm} + K_{\soqdm}} \, \frac{10^{-3}\zeta^\l C_{\hdl}}{10^{-3}\zeta^\l C_{\hdl} + K_{\hdl}}, $$
	with
	$v_\text{max} = 1.45\times 10^{-5}$ $\mathrm{mol \cdot L^{-1} \cdot \text{year}^{-1}}$,
	$K_{\soqdm} = 10^{-4}$ $\mathrm{mol \cdot L^{-1}}$
	and $K_{\hdl} = \varepsilon$ $\mathrm{mol \cdot L^{-1}}$.
	
	Table~\ref{tab:hd-init} in Appendix~\ref{sec:values-h2} presents the initial conditions of the system, which is assumed to be at equilibrium. In Table~\ref{tab:hd-const} in Appendix~\ref{sec:values-h2}, we give the equilibrium constants linked to the chemical reactions. The remaining constants are chosen to be identical to those used in the $\cod$ test case.

	\subsubsection{Vertical well model}
	
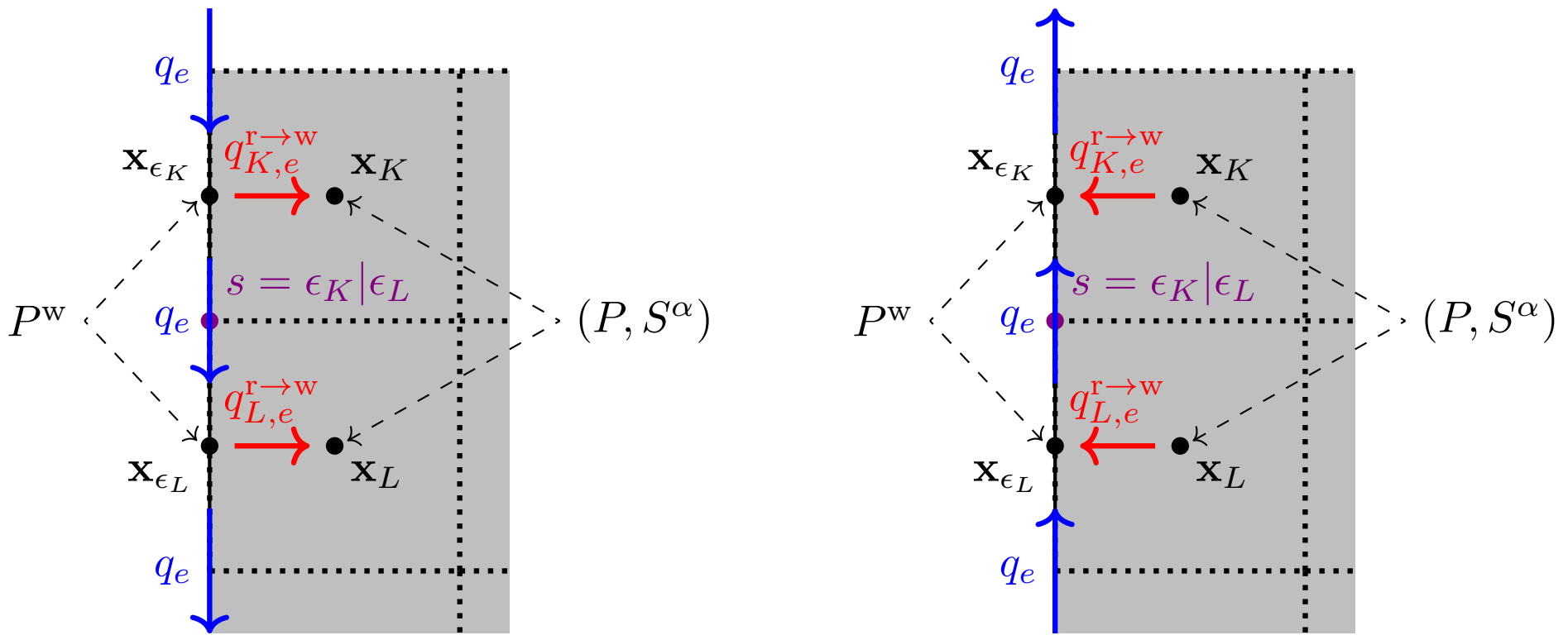
\begin{figure}[H]
	\centering
	\subfloat[Injection]{
		\begin{tikzpicture}[scale=2]
			\definecolor{dom2}{RGB}{191, 191, 191}      
			
			\draw[fill=dom2, draw=none] (0, -0.25) rectangle (1.2, 2);
			\draw[step=1cm, very thick, dotted] (0,0) grid (1,2);
			\draw[very thick, dotted] (1,0) -- (1,-0.25);
			\draw[very thick, dotted] (1,0) -- (1+0.2,0);
			\draw[very thick, dotted] (1,1) -- (1+0.2,1);
			\draw[very thick, dotted] (1,2) -- (1+0.2,2);
			
			\draw[thick] (0,0) -- (0,2);
			
			\fill (0,0.5) circle (1pt) node[below left] {$\mathbf{x}_{\epsilon_L}$};
			\fill (0,1.5) circle (1pt) node[above left] {$\mathbf{x}_{\epsilon_K}$};
			\fill[violet] (0,1) circle (1pt) node[violet, above right] {$s=\epsilon_K|\epsilon_L$};
			
			\fill (0.5,0.5) circle (1pt) node[below right] {$\mathbf{x}_L$};
			\fill (0.5,1.5) circle (1pt) node[above right] {$\mathbf{x}_K$};
			
			\draw[->, very thick, red] (0.1,0.5) -- (0.5-0.1,0.5) node[midway,above] {$q_{L,e}^{\mathrm{r}\to\mathrm{w}}$};
			\draw[->, very thick, red] (0.1,1.5) -- (0.5-0.1,1.5) node[midway,above] {$q_{K,e}^{\mathrm{r}\to\mathrm{w}}$};
			
			\draw[<-, very thick, blue] (0,-0.25) -- (0,0.25) node[midway,left] {$q_{e}$};
			\draw[<-, very thick, blue] (0,0.75) -- (0,1.25) node[midway,left] {$q_{e}$};
			\draw[<-, very thick, blue] (0,1.75) -- (0,2.25) node[midway,left] {$q_{e}$};
			
			\draw[<-, dashed] (0.5+0.05,0.5+0.02) -- (1.4,1) node[right] {$(P,S^\alpha)$};
			\draw[<-, dashed] (0.5+0.05,1.5-0.02) -- (1.4,1);
			
			\draw[<-, dashed] (0-0.05,0.5+0.02) -- (-0.5,1) node[left] {$P^\mathrm{w}$};
			\draw[<-, dashed] (0-0.05,1.5-0.02) -- (-0.5,1);
		\end{tikzpicture}
	}
	\qquad
	\subfloat[Production]{
		\begin{tikzpicture}[scale=2]
			\definecolor{dom2}{RGB}{191, 191, 191}      
			
			\draw[fill=dom2, draw=none] (0, -0.25) rectangle (1.2, 2);
			\draw[step=1cm, very thick, dotted] (0,0) grid (1,2);
			\draw[very thick, dotted] (1,0) -- (1,-0.25);
			\draw[very thick, dotted] (1,0) -- (1+0.2,0);
			\draw[very thick, dotted] (1,1) -- (1+0.2,1);
			\draw[very thick, dotted] (1,2) -- (1+0.2,2);
			
			\draw[thick] (0,0) -- (0,2);
			
			\fill (0,0.5) circle (1pt) node[below left] {$\mathbf{x}_{\epsilon_L}$};
			\fill (0,1.5) circle (1pt) node[above left] {$\mathbf{x}_{\epsilon_K}$};
			\fill[violet] (0,1) circle (1pt) node[violet, above right] {$s=\epsilon_K|\epsilon_L$};
			
			\fill (0.5,0.5) circle (1pt) node[below right] {$\mathbf{x}_L$};
			\fill (0.5,1.5) circle (1pt) node[above right] {$\mathbf{x}_K$};
			
			\draw[<-, very thick, red] (0.1,0.5) -- (0.5-0.1,0.5) node[midway,above] {$q_{L,e}^{\mathrm{r}\to\mathrm{w}}$};
			\draw[<-, very thick, red] (0.1,1.5) -- (0.5-0.1,1.5) node[midway,above] {$q_{K,e}^{\mathrm{r}\to\mathrm{w}}$};
			
			\draw[->, very thick, blue] (0,-0.25) -- (0,0.25) node[midway,left] {$q_{e}$};
			\draw[->, very thick, blue] (0,0.75) -- (0,1.25) node[midway,left] {$q_{e}$};
			\draw[->, very thick, blue] (0,1.75) -- (0,2.25) node[midway,left] {$q_{e}$};
			
			\draw[<-, dashed] (0.5+0.05,0.5+0.02) -- (1.4,1) node[right] {$(P,S^\alpha)$};
			\draw[<-, dashed] (0.5+0.05,1.5-0.02) -- (1.4,1);
			
			\draw[<-, dashed] (0-0.05,0.5+0.02) -- (-0.5,1) node[left] {$P^\mathrm{w}$};
			\draw[<-, dashed] (0-0.05,1.5-0.02) -- (-0.5,1);
		\end{tikzpicture}
	}
	\caption{Discretization of the well at the left boundary $r=r_{\mathrm{w}}$ of the radial cross-section.}
	\label{fig:well-discretization}
	\end{figure}

	In this section, we present a simplified model of a vertical well that allows for cycles of chemical species injection into and withdrawal from the reservoir.
	Since the characteristic flow time scale in the well is much shorter than that in the reservoir, we assume that the flow in the well is at steady state and satisfies the molar conservation of each species along the well \cite{beaude_parallel_2018,armandine_les_landes_two-phase_2023}.
	It is also assumed that crossflow between injection and production intervals is excluded in the sense that all well fluxes are either positive for production or negative for injection.
	In radial geometry, the well is discretized using a set of edges and we denote by $\mathcal{M}^\mathrm{w}$ the set of reservoir cells connected to it (see Figure~\ref{fig:well-discretization}). For a given cell $K$ in $\mathcal{M}^\mathrm{w}$, we let $\mathbf{x}_{\epsilon_K}$ be the center of the well edge intersecting this cell.
	
	Using an upwind approximation of the mobility, the species molar flow rate between the reservoir and the well is expressed	as follows
	\begin{equation}
	q_{K,e}^{\text{r}\to \text{w}} := \begin{cases}
		\zeta^{g}(P^\text{w}_K) \, \eta^{g}(1)\, \text{WI}\, \min\big(0, P_K-P^\text{w}_K\big) & \text{if}\ e=\hdg, \\
		0 & \text{otherwise},
	\end{cases} \leq 0.
	\label{eq:qrwk-injection}
	\end{equation}
	during pure $\hdg$ injection and 
	\begin{equation}
	q_{K,e}^{\text{r}\to \text{w}} := \zeta^{\alpha(e)}\big(P_K,C_K^{\alpha(e)}\big)\, \eta^{\alpha(e)}\big(S^{\alpha(e)}_K\big)C_{e,K}\, \text{WI}\, \max\big(0, P_K- P^\text{w}_K\big) \geq 0,
	\label{eq:qrwk-production}
	\end{equation}
	during production, where $q_{K,e}^{\text{r}\to \text{w}}$ denotes the molar flow rate of the species $e$ between the reservoir and the well [$\mathrm{mol}\cdot\mathrm{s}^{-1} $] and $\text{WI}$ is the well index [$\mathrm{m}$].
	We take a reference pressure in a given cell of $\mathcal{M}^\mathrm{w}$ (bottom-hole pressure) $P^\mathrm{w}_\mathrm{bh}$ and express the well pressures, assuming an hydrostatic pressure drop, by
	\begin{equation}
	P^\mathrm{w}_K = P^\mathrm{w}_\mathrm{bh} + \Delta P^\mathrm{w}_K.
	\label{eq:well-Pbh}
	\end{equation}
	We impose either a prescribed total flow rate or a prescribed bottom-hole pressure through the following complementary constraints: during injection ($\overline{Q}\leq 0$)
	\begin{equation}
	\min\big( Q_h^{\text{r}\to \text{w}}-\overline{Q}, \overline{P} - P^\mathrm{w}_\mathrm{bh} \big) = 0,
	\label{eq:well-inj}
	\end{equation}
	and during production ($\overline{Q}\geq 0$)
	\begin{equation}
	\min\big( \overline{Q}-Q_h^{\text{r}\to \text{w}}, P^\mathrm{w}_\mathrm{bh} - \overline{P} \big) = 0,
	\label{eq:well-prod}
	\end{equation}
	with $Q_h^{\text{r}\to \text{w}}$ the discretized total molar flow rate [$\mathrm{mol}\cdot\mathrm{s}^{-1}$], $\overline{Q}$ the target total molar flow rate [$\mathrm{mol}\cdot\mathrm{s}^{-1}$], and $\overline{P}$ the maximum bottom-hole pressure for injection and minimum bottom-hole pressure for production [$\mathrm{Pa}$]. We have
	$$ Q_h^{\text{r}\to \text{w}} := \sum_{K\in\mathcal{M}^\mathrm{w}} \sum_{e^\text{pr}\in\E^\text{pr}\cap\E^\text{f}} \sum_{e\in\E^\text{f}} E_{e^\text{pr},e}\,q_{K,e}^{\text{r}\to \text{w}}. $$
	It remains to determine $\Delta P^\mathrm{w}_K$ in \eqref{eq:well-Pbh}. To this end, we assume the bottom-hole pressure $P^{\text{w},n-1}_{\mathrm{bh}}$ and the mixture mass density in the well $\rho^{\mathrm{w},n-1}$ to be known at iteration $n-1$, and we suppose hydrostatic pressure within the well. Assuming that $\rho^{\mathrm{w},n-1}=\rho^{\mathrm{w},n-1}_K$ is constant on each cell $K$, we have
	$$ \Delta P^{\text{w},n-1}_K = - \sum_{s=\epsilon_L|\epsilon_M\in\mathcal{V}_\mathrm{bh}^\mathrm{w}(K)} g\, \frac{\rho^{\mathrm{w},n-1}_L+\rho^{\mathrm{w},n-1}_M}{2} \,(z_{\epsilon_L}-z_{\epsilon_M}), $$
	with $z_{\epsilon_L}$ the $z$-coordinate of the center of the edge $\epsilon_L$ and $\mathcal{V}_\mathrm{bh}^\mathrm{w}(K)$ is the set of well's vertices between the bottom-hole edge and $\epsilon_K$.
	In this way, it is sufficient to know the bottom-hole pressure and the mass density in each cell to compute the pressure in the well $P^\text{w}_K$. The well equation to be solved is: for given $\Delta P^{\text{w},n-1}_K$, find $P^{\text{w}}_{\mathrm{bh}}$ such that, 
	during injection
	\begin{equation*}
	\min\Bigg( \sum_{e^\text{pr}\in\E^\text{pr}\cap\E^\text{f}} \sum_{K\in\mathcal{M}^\mathrm{w}} \sum_{e\in\E^\text{f}} E_{e^\text{pr},e}\, q_{K,e}^{\text{r}\to \text{w}}\big(P^{\text{w}}_{\mathrm{bh}} + \Delta P^{\text{w},n-1}_K\big)-\overline{Q}, \overline{P} - P^\mathrm{w}_\mathrm{bh} \Bigg) = 0,
	\end{equation*}
	and during production
	\begin{equation*}
	\min\Bigg( \overline{Q}-\sum_{e^\text{pr}\in\E^\text{pr}\cap\E^\text{f}} \sum_{K\in\mathcal{M}^\mathrm{w}} \sum_{e\in\E^\text{f}} E_{e^\text{pr},e}\, q_{K,e}^{\text{r}\to \text{w}}\big(P^{\text{w}}_{\mathrm{bh}} + \Delta P^{\text{w},n-1}_K\big), P^\mathrm{w}_\mathrm{bh} - \overline{P} \Bigg) = 0.
	\end{equation*}
	This equation is added to the global implicit system.
	From an implementation point of view, we introduce an additional algebraic degree of freedom associated with the well equation in the Jacobian, with the bottom-hole pressure of the well corresponding to the pressure column of this degree of freedom. In this way, we preserve the performance of the CPR-AMG solver by preconditioning the Jacobian matrix with the inverse of the pressure block.
	
	During injection, we know what is being injected, and then the density in the well is based on the previous time step well pressures. In our case, we have
	$$ \rho_K^{\mathrm{w},n-1} = \rho^g\big(P^{\text{w},n-1}_K\big), $$
	where $P^{\text{w},n-1}_K$ is the pressure in the well at the previous iteration. During production,
	we use the mixture density in the cell as a first approximation of the well mixture density, and we have
	$$ \rho_K^{\mathrm{w},n-1} = S^{\l,n-1}_K \rho^\l\big(P^{n-1}_K\big) + S^{g,n-1}_K \rho^g\big(P^{n-1}_K\big). $$
	
	\subsubsection{Well index computation in radial geometry}
	
\begin{figure}[H]
	\centering
	\begin{tikzpicture}[scale=1]
		\definecolor{dom2}{RGB}{191, 191, 191}      
		
		\def\R{3}      
		\def\r{1}      
		\def\ry{0.6}     
		\def\h{3.5}      
		
		\draw[fill=dom2, draw=none] (\r, 0) rectangle (\R, \h);
		
		
		\draw (-\R,0) -- (-\R,\h);
		\draw[very thick, dotted]  (\R,0) -- (\R,\h);
		
		\draw (0,\h) ellipse (\R cm and \ry cm);
		
		\draw (0,0) ellipse (\R cm and \ry cm);

		
		\draw (-\r,0) -- (-\r,\h);
		\draw[very thick, dotted] (\r,0) -- (\r,\h);
		
		\draw (0,\h) ellipse (\r cm and \ry*0.3 cm);
		
		\draw (0,0) ellipse (\r cm and \ry*0.3 cm);
		
		\draw[->] (0,\h*0.75) -- (\r,\h*0.75) node[midway, below] {$r_\mathrm{w}$};
		
		\draw[->] (0,\h*0.25) -- (\R*0.5+\r*0.5,\h*0.25) node[midway, above] {$r_K$};
		
		\draw[dashed] (\R*0.5+\r*0.5,\h*0.5) -- (\R*0.5+\r*0.5,\h*0.25);
		
		\draw[<->] (\R+0.5,0) -- (\R+0.5,\h) node[midway, right] {$h$};
		
		\draw (0,0) -- (0,\h);
		
		\fill (\r,\h*0.5) circle (2pt) node[left] {$\mathbf{x}_{\epsilon_K}$};
		\fill (\R*0.5+\r*0.5,\h*0.5) circle (2pt) node[above] {$\mathbf{x}_K$};
		
		\draw[very thick, dotted] (\r,0) -- (\R,0);
		\draw[very thick, dotted] (\r,\h) -- (\R,\h);
		
	\end{tikzpicture}		
	\caption{Discretization around the well in radial geometry.}
	\label{fig:well-radial-discretization}
	\end{figure}
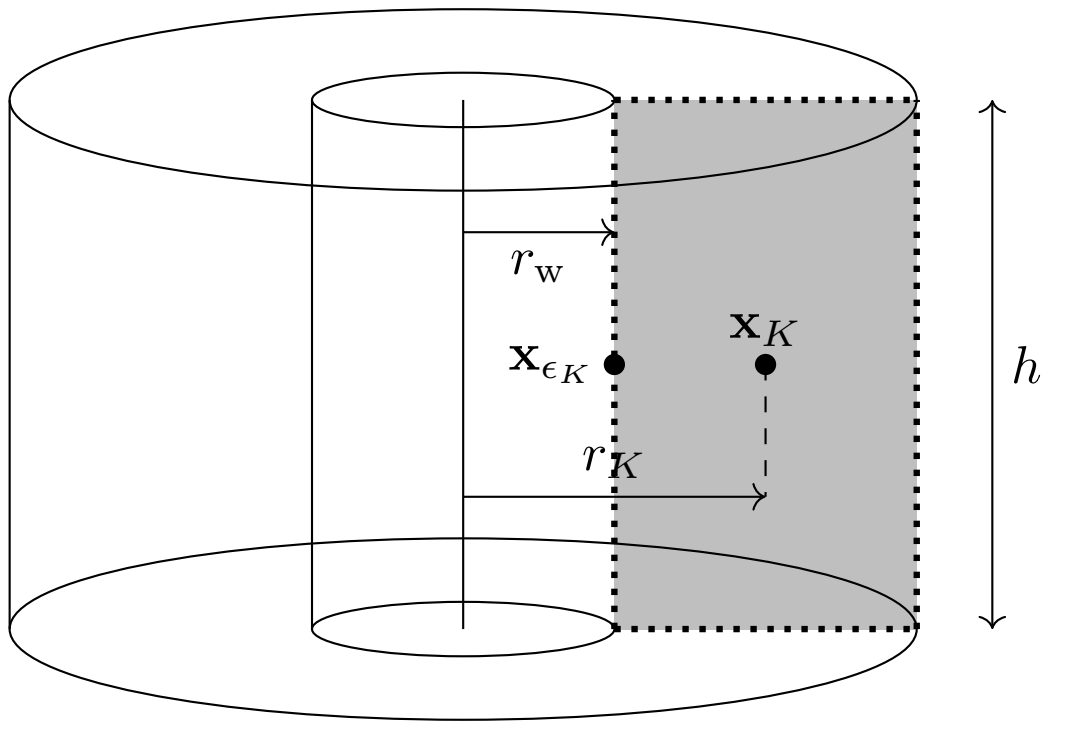

	To compute the total molar flow rate entering the well, denoted by $Q_h^{\mathrm{r}\to\mathrm{w}}$ in \eqref{eq:well-inj} and \eqref{eq:well-prod}, it remains to determine the well index $\mathrm{WI}$ corresponding to the considered geometry and numerical scheme. We first introduce the analytical reference solution, assuming a homogeneous and isotropic medium (i.e. $\mathbb{K} = k\,\mathrm{Id}$), for a vertical well in a radially unbounded reservoir of thickness $h$, with a single chemical species and a single phase
	\cite{peaceman_interpretation_1978,chen_well_2009}. The radial pressure solution is given by
	\begin{equation}
	P(r) = P(r_\mathrm{w}) + \frac{1}{2\pi h}\frac{1}{\zeta}\frac{\mu}{k}Q^{\mathrm{r}\to\mathrm{w}}\ln\left(\frac{r}{r_\mathrm{w}}\right),\quad r>0
	\label{eq:radial-pressure-solution}
	\end{equation}
	which satisfies the molar conservation equation
	\begin{equation}
	-\operatorname{div}\bigg(\zeta\frac{k}{\mu}\nabla P\bigg) = Q^{\mathrm{r}\to\mathrm{w}}\delta_\mathrm{w},
	\label{eq:radial-conservation}
	\end{equation}
	where $r_w$ is the well radius, $\zeta$ is the molar density [$\mathrm{mol\cdot m^{-3}}$], $\mu$ is the viscosity $[\mathrm{Pa}\cdot\mathrm{s}]$, $Q^{\mathrm{r}\to\mathrm{w}}$ denotes the total molar flow rate entering the well $[\mathrm{mol}\cdot\mathrm{s}^{-1}]$, and $\delta_w$ denotes the Dirac measure concentrated at the well location, i.e., at the origin of the radial domain. In such settings, note that $Q_h^{\mathrm{r}\to\mathrm{w}}$ in \eqref{eq:well-inj}--\eqref{eq:well-prod} and $q_{K,e}^{\text{r}\to \text{w}}$ in \eqref{eq:qrwk-injection}--\eqref{eq:qrwk-production} are equal to $Q^{\mathrm{r}\to\mathrm{w}}$. The well index is then determined by requiring the discrete well-reservoir coupling to reproduce the pressure-flow-rate relationship of the analytical radial solution (see Figure~\ref{fig:well-radial-discretization}).
	
	The flow rate between the reservoir and the well \eqref{eq:qrwk-injection}--\eqref{eq:qrwk-production} reads
	$$ \zeta\frac{k}{\mu}\mathrm{WI}\,\left(P_K-P^\mathrm{w}_K\right) = Q^{\mathrm{r}\to\mathrm{w}}. $$
	Using the analytical solution given in \eqref{eq:radial-pressure-solution}, we have $P^\mathrm{w}_K = P(r_\mathrm{w})$ and
	$$ P_K = P(r_K) = P^\mathrm{w}_K + \frac{1}{2\pi h}\frac{1}{\zeta}\frac{\mu}{k}Q^{\mathrm{r}\to\mathrm{w}} \ln\left(\frac{r_K}{r_\mathrm{w}}\right). $$
	Replacing these expressions in the discretization computed above leads to
	$$ \mathrm{WI} = \frac{2\pi h}{\ln({r_K}/{r_\mathrm{w}})}. $$
	
	\subsubsection{Results}
	
\begin{figure}[H]
	\centering
	\subfloat{\includegraphics[width=7.5cm]{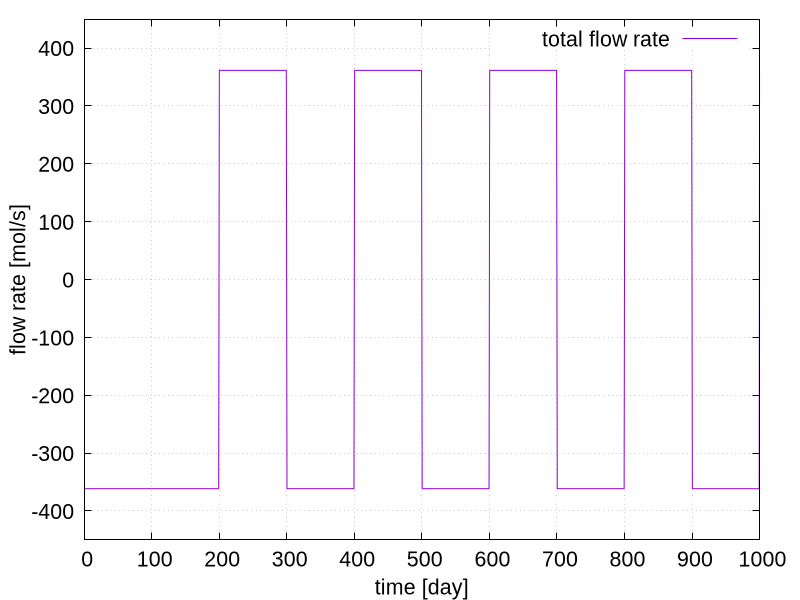}}
	\subfloat{\includegraphics[width=7.5cm]{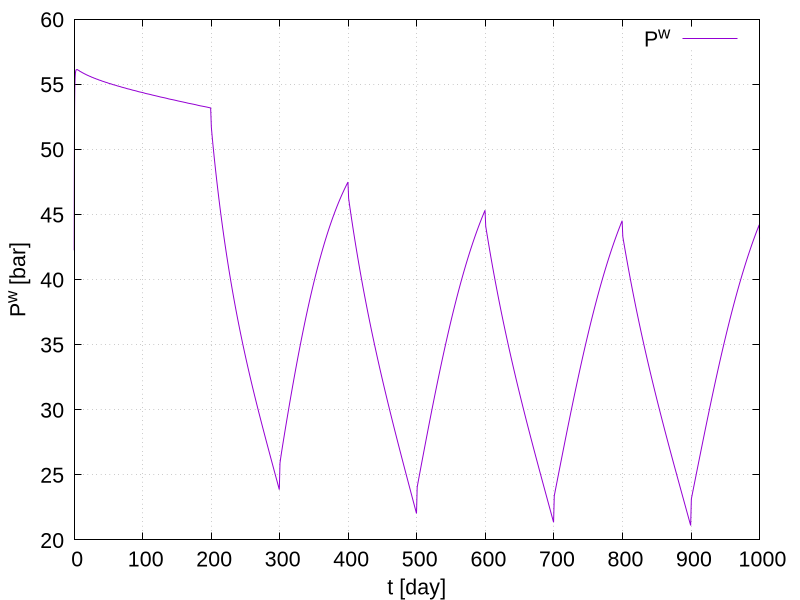}}
	\quad
	\subfloat{\includegraphics[width=7.5cm]{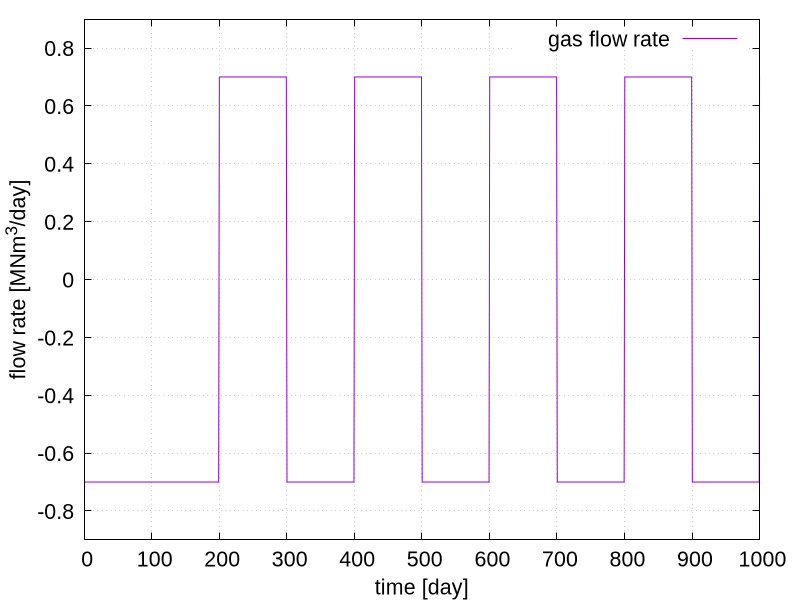}}
	\subfloat{\includegraphics[width=7.5cm]{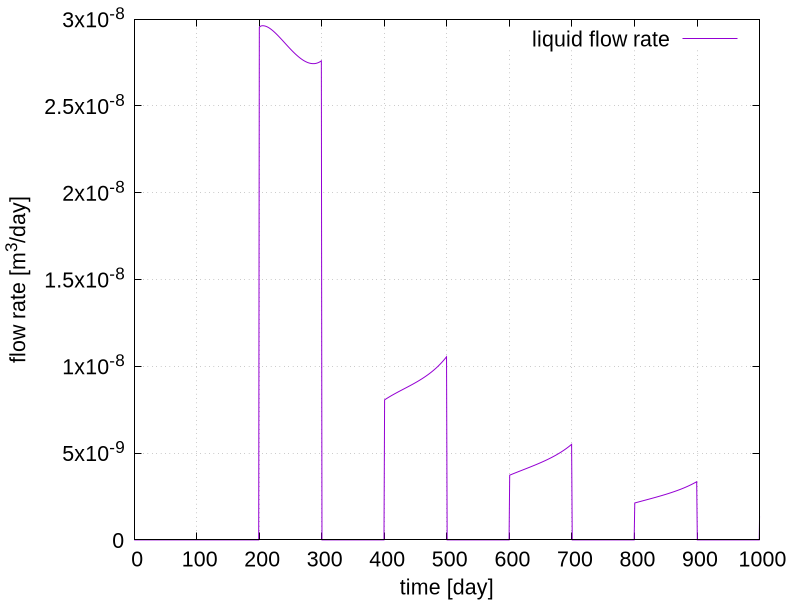}}
	
	\caption{$\hd$ storage test case. Top-left: total molar flow rate between the reservoir and the well $Q_h^{\text{r}\to \text{w}}$ over time. Top-right: well bottom hole pressure $P^{\text{w}}_{\mathrm{bh}}$ over time. Bottom-left: gas flow rate between the reservoir and the well over time. Bottom-right: liquid flow rate between the reservoir and the well over time.}
	\label{fig:pwk_wrt_time}
	\end{figure}

	The numerical simulation is performed in radial geometry on a non-uniform mesh containing $4096$ elements, with a refined discretization near the well. The simulation consists of alternating 100-day cycles of pure $\hdg$ injection and production, with the total molar flow rate prescribed during each cycle, see Figure~\ref{fig:pwk_wrt_time}. The initial injection begins at the start of the simulation and lasts 200 days, thereby establishing a cushion gas inventory for subsequent cycles and preventing liquid production.
	
	As shown in Figure~\ref{fig:pwk_wrt_time}, the bottom-hole pressure required during injection increases as the reservoir pressure equilibrates with the pressure in the well. The same effect is observed during production, causing the bottom-hole pressure to decrease to maintain the prescribed total molar flow rate. Finally, a small amount of liquid is produced near the end of withdrawal cycles. This behavior may be influenced by the absence of water vaporization in the chemical model.

	Table~\ref{tab:benchmark-hd} gives the total number of time steps performed by the code, the average number of nonlinear iterations per time step, the average number of linear iterations per Newton iteration, and the {CPU time}. No disappearance of the aqueous phase is observed, which simplifies the simulation and reduces the occurrence of time-step failures. In addition, no flash calculations are required, which speeds up the computation despite the presence of more chemical species than in the $\cod$ test case.
	
	\begin{table}[H]
	\centering
	\begin{tabular}{|c|c|c|c|c|c|}
		\hline
		\rowcolor{gray!40} TS (failed) & NI & LI / NI & F (failed) & CPU & Architecture \\
		\hline
		1077 (20) & 10.50 & 13.04 & 0 (0) & 12052 & Intel Core Ultra 7 255H @2.0 GHz\\
		\hline
	\end{tabular}\vspace{0.2cm}		
	\caption{Numerical performance of our method for $\hd$ storage application. See caption of Table~\ref{tab:benchmark-cod} for the meaning of the abbreviations.}
	\label{tab:benchmark-hd}
	\end{table}

\begin{figure}[H]
	\centering
	\subfloat{\includegraphics[height=5.5cm]{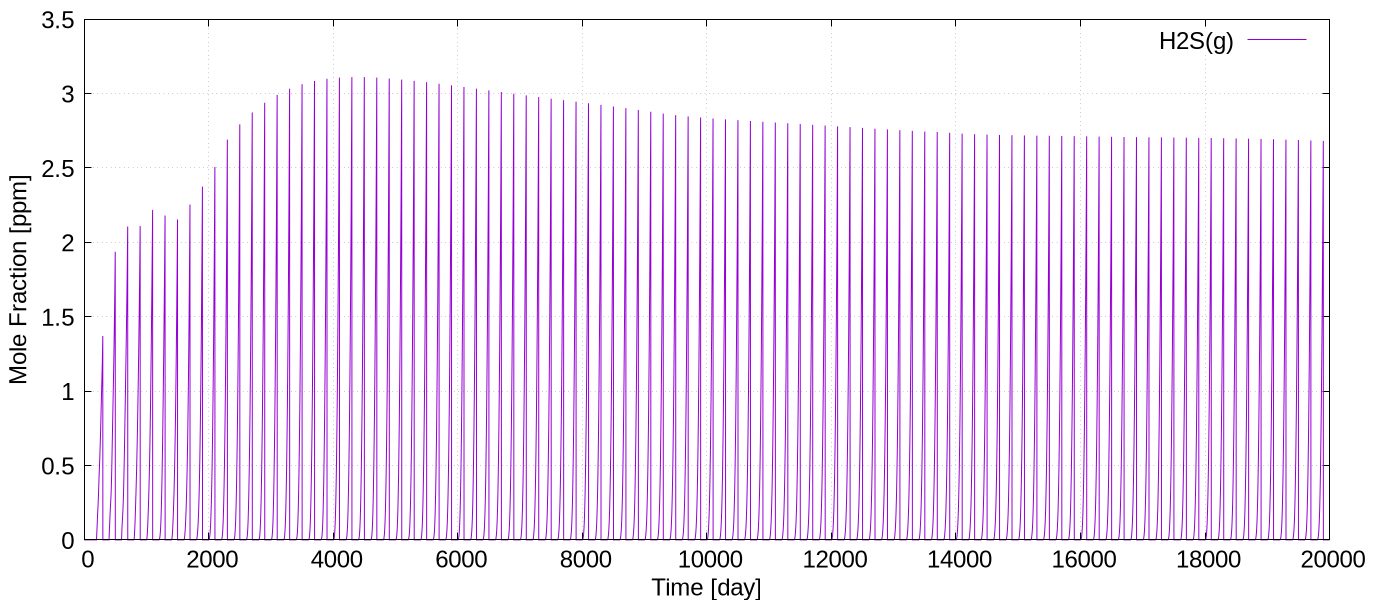}}
	\quad
	\subfloat{\includegraphics[height=5.5cm]{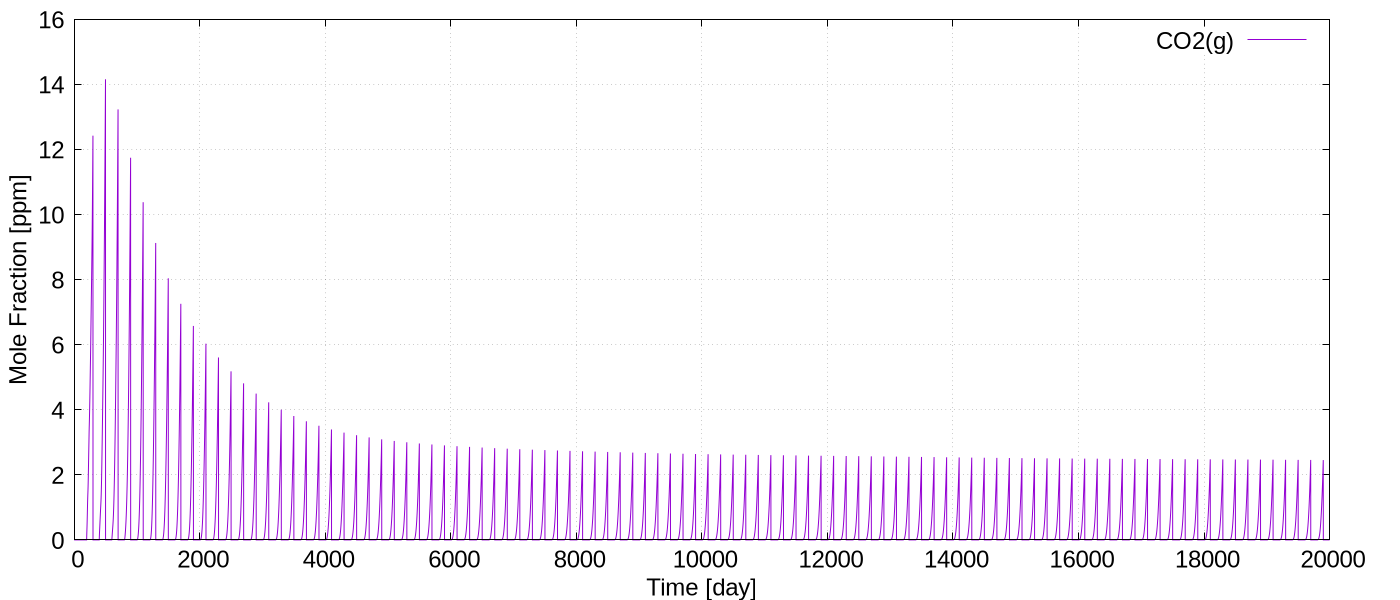}}
	
	\caption{$\hd$ storage test case: Production of $\hdsg$ and of $\codg$ over days.}
	\label{fig:h2sg-production-ppm}
	\end{figure}

	Figure~\ref{fig:h2sg-production-ppm} shows the produced mole fractions of $\hdsg$ and of $\codg$ over time for a period of 20,000 days. The produced $\hdsg$ mole fraction initially increases over successive injection–withdrawal cycles, reaches a maximum of a few  ppm, and subsequently decreases toward an approximately periodic regime. Conversely, the produced $\codg$ mole fraction decreases markedly during the early cycles (except for the two  first withdrawals) before approaching a nearly periodic response.
	
\begin{figure}[H]
	\centering
	\includegraphics[height=6.8cm]{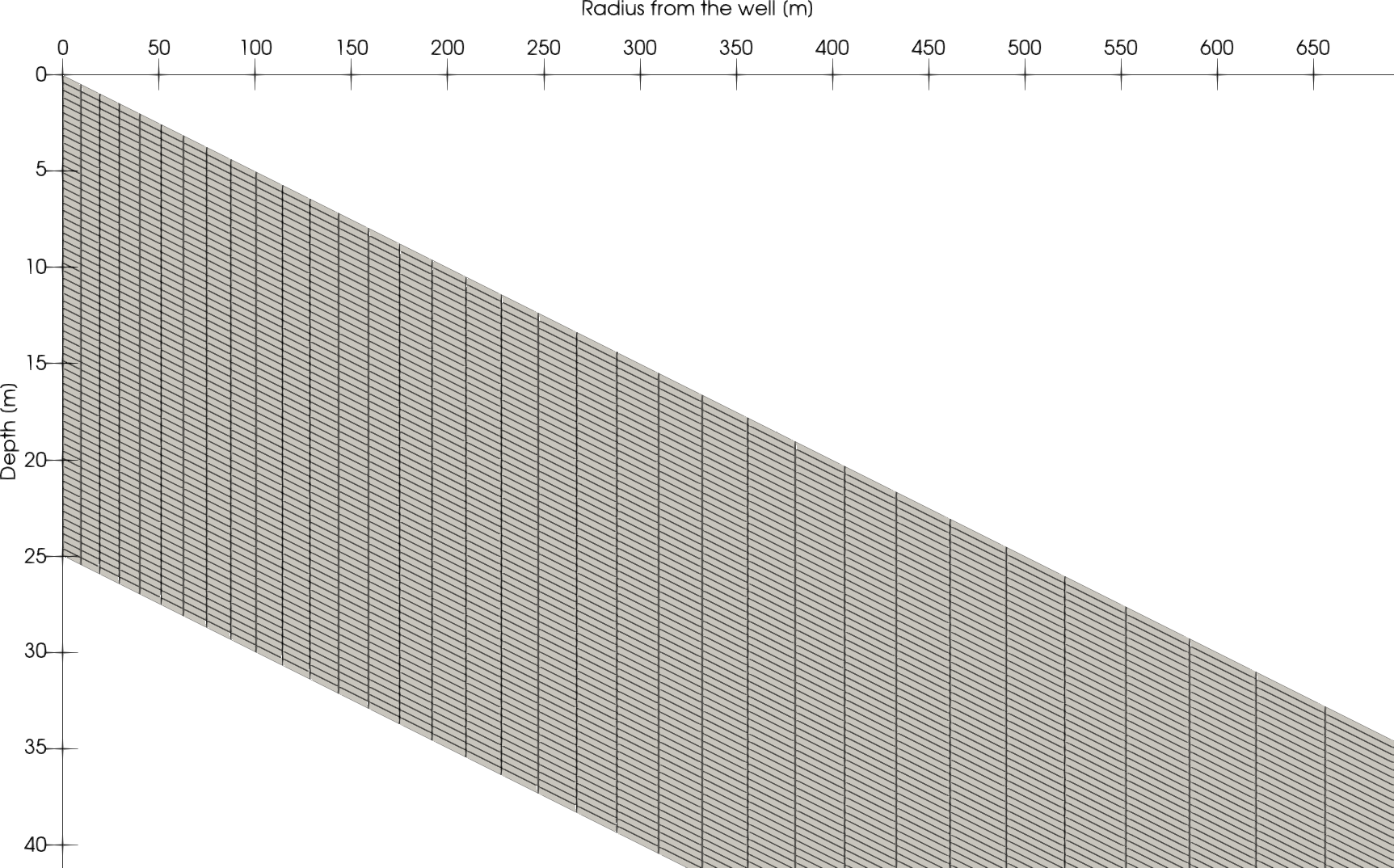}
	
	\caption{$\hd$ storage test case: Close-up of the radial cross-section mesh, highlighting the local mesh refinement around the well boundary.}
	\label{fig:h2-mesh}
	\end{figure}

\begin{figure}[H]
	\centering
	\subfloat[molar concentration of $\hdg$ at $t=200$ days]{
		\includegraphics[width=\2dimgw]{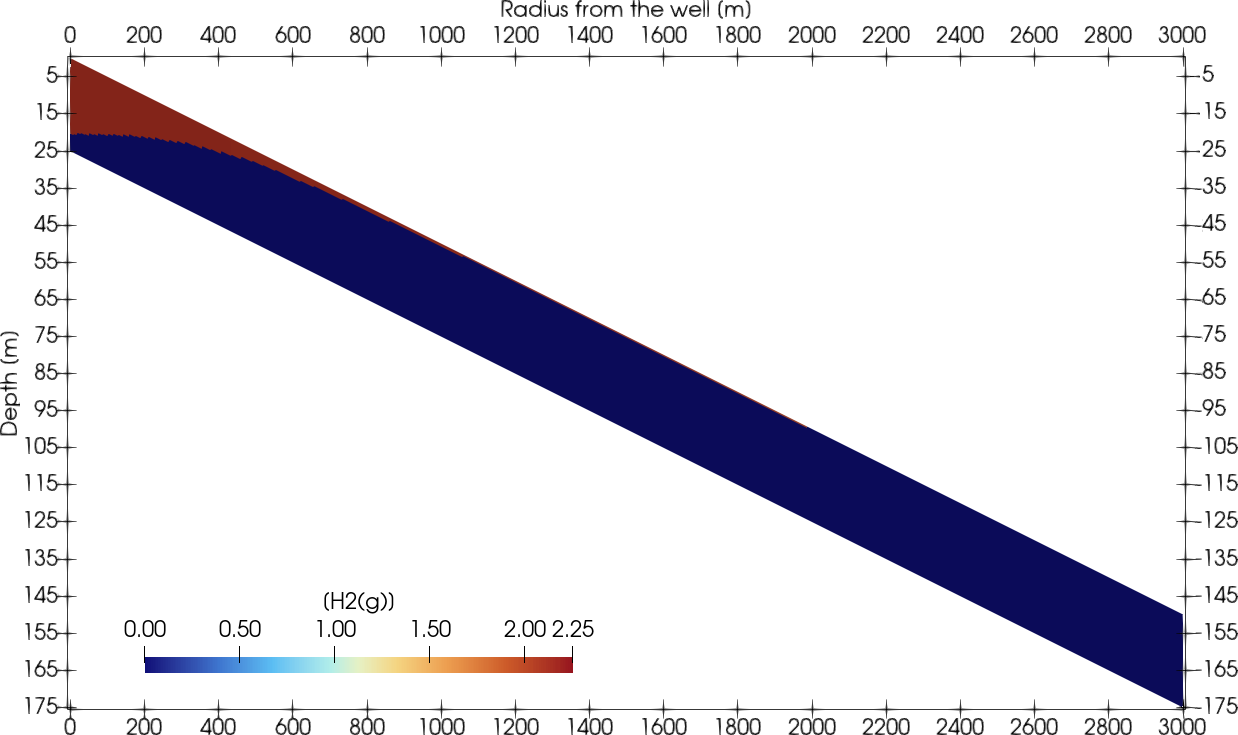}
	}
	\subfloat[molar concentration of $\hdg$ at $t=300$ days]{
		\includegraphics[width=\2dimgw]{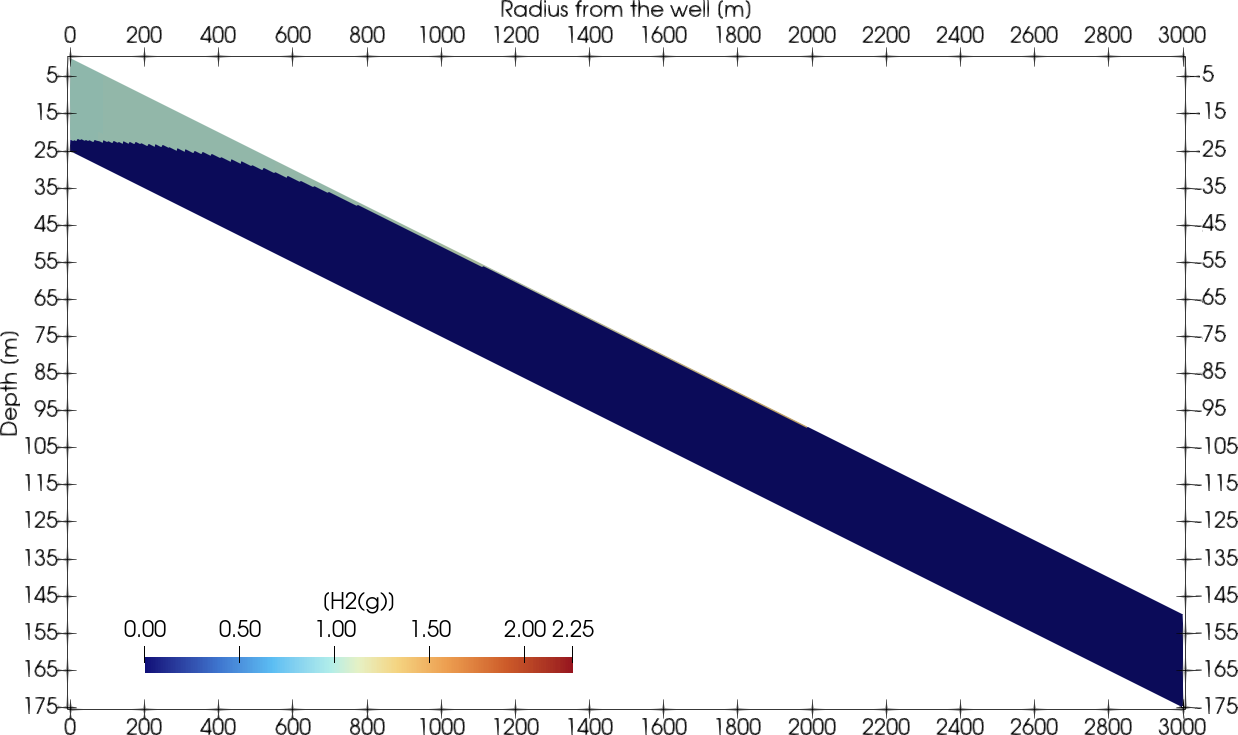}
	}
	\quad
	\subfloat[molar concentration of $\hdsg$ at $t=200$ days]{
		\includegraphics[width=\2dimgw]{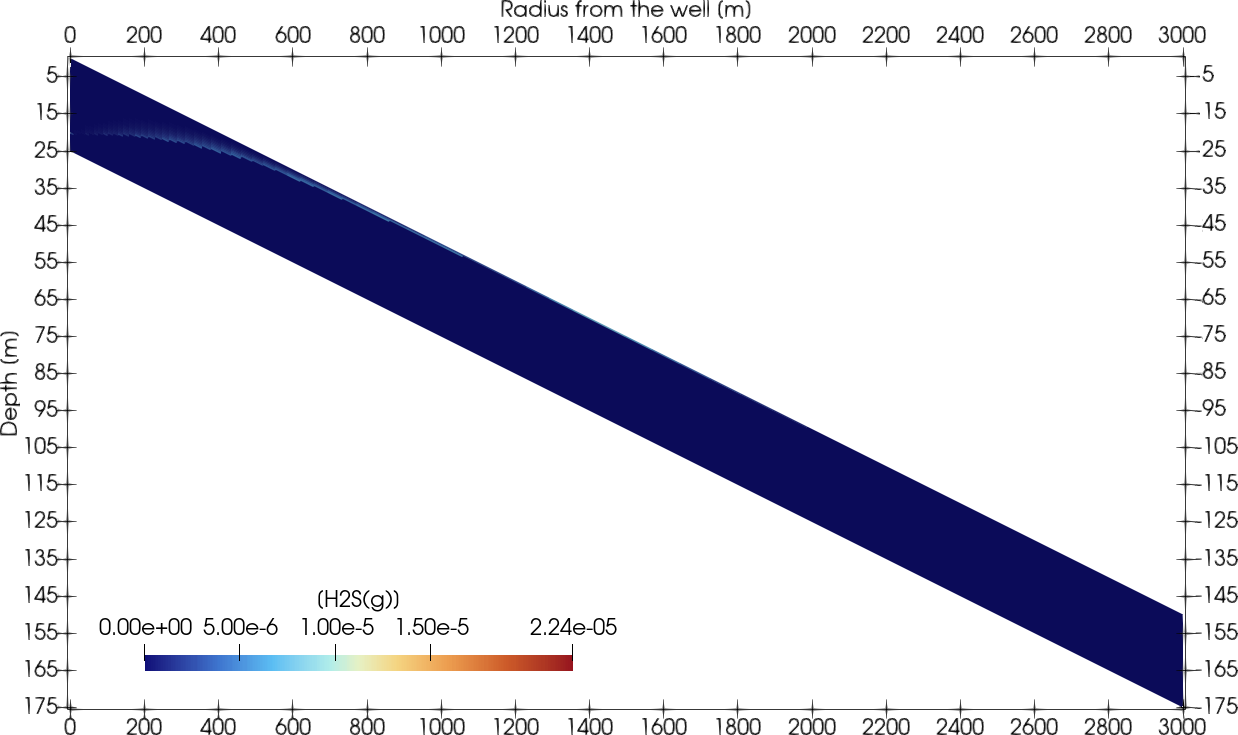}
	}
	\subfloat[molar concentration of $\hdsg$ at $t=300$ days]{
		\includegraphics[width=\2dimgw]{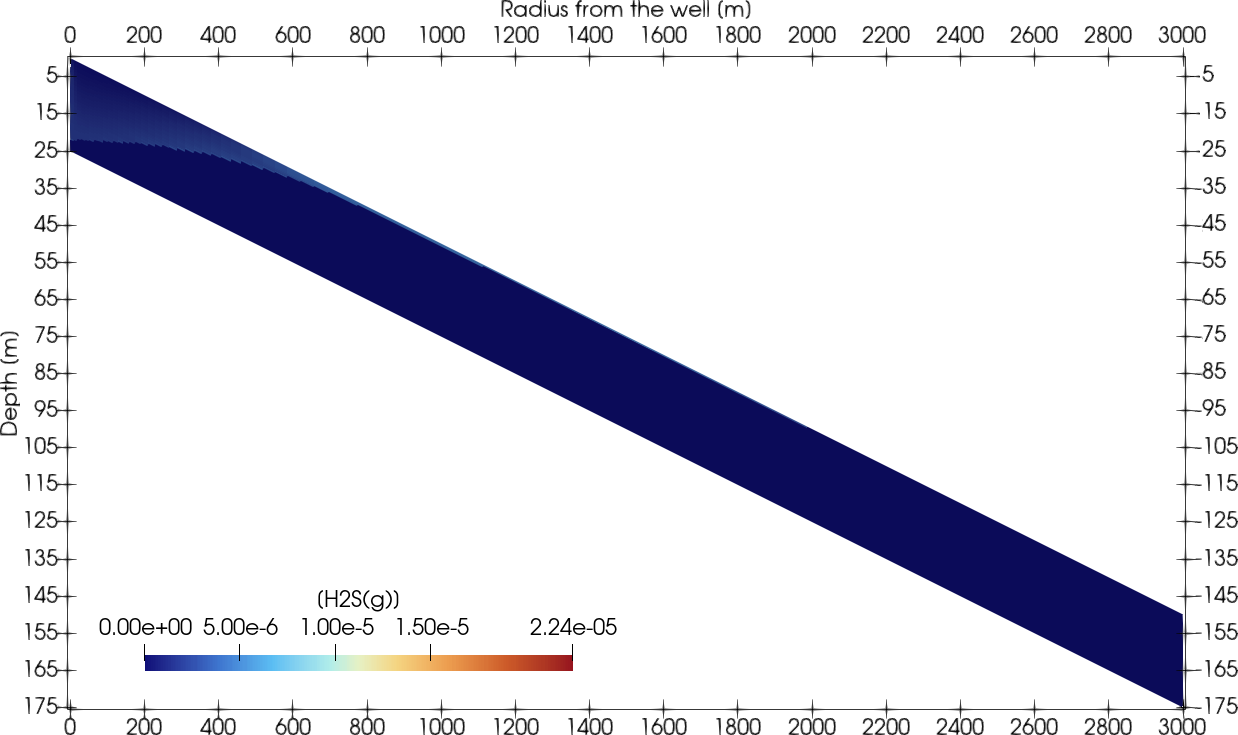}
	}
	\quad
	\subfloat[pH at $t=200$ days]{
		\includegraphics[width=\2dimgw]{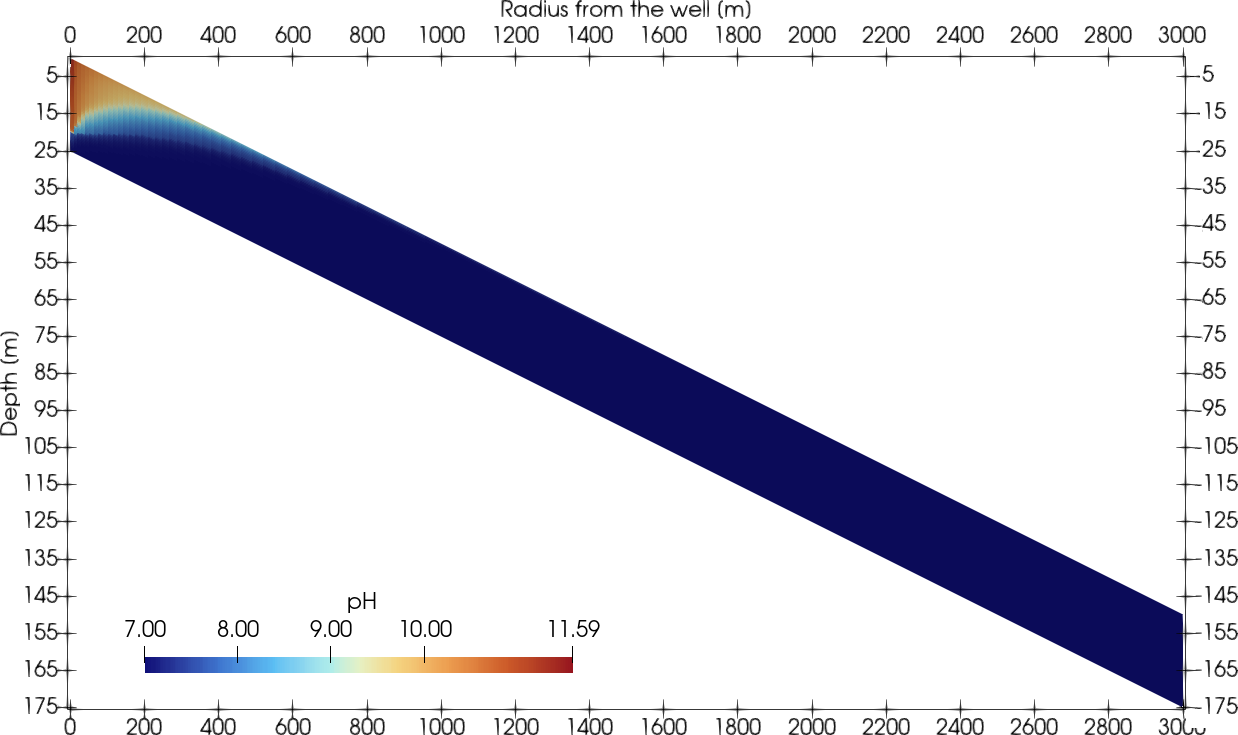}
	}
	\subfloat[pH at $t=300$ days]{
		\includegraphics[width=\2dimgw]{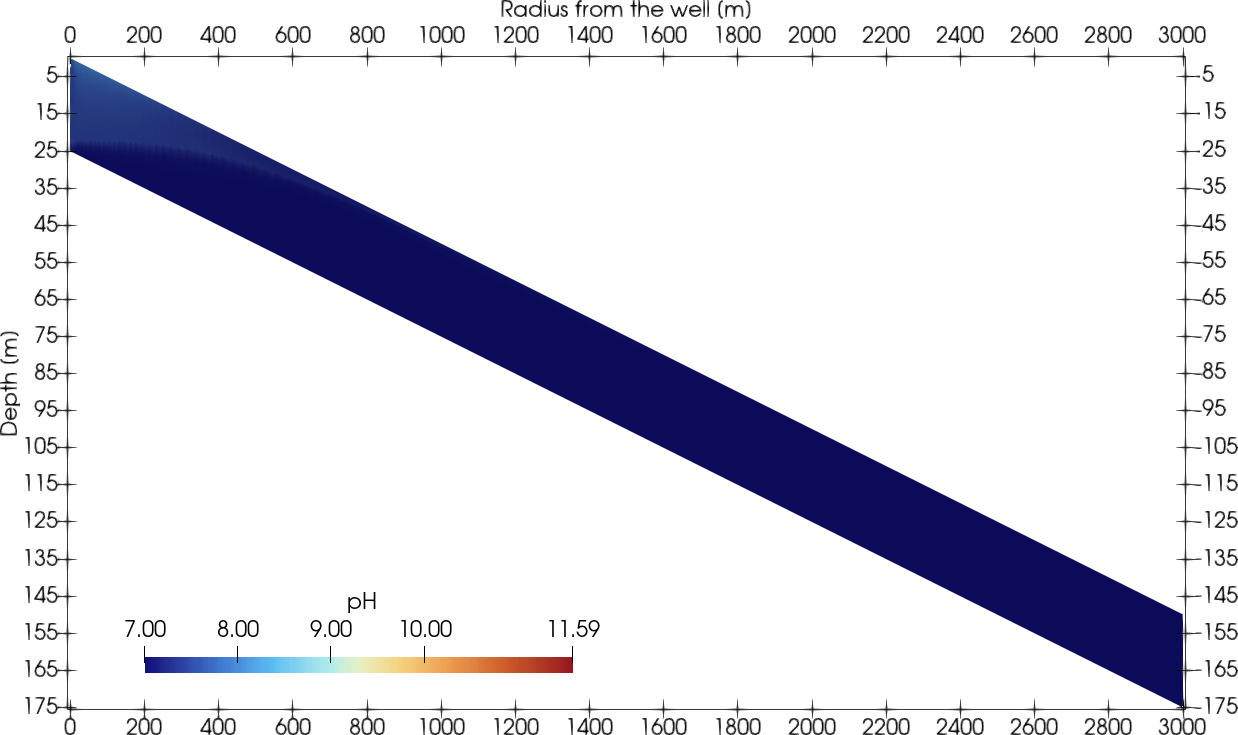}
	}
	\caption{Results of the gaseous $\hd$ injection and withdrawal simulation in a reservoir with radial geometry at times $t = 200$ days and $t = 300$ days. The molar concentrations of $\hdg$ and $\hdsg$, as well as pH are shown.}
	\label{fig:ben-hd-rad-first}
	\end{figure}

\begin{figure}[H]
	\centering
	\subfloat[molar concentration of $\hdg$ at $t=900$ days]{
		\includegraphics[width=\2dimgw]{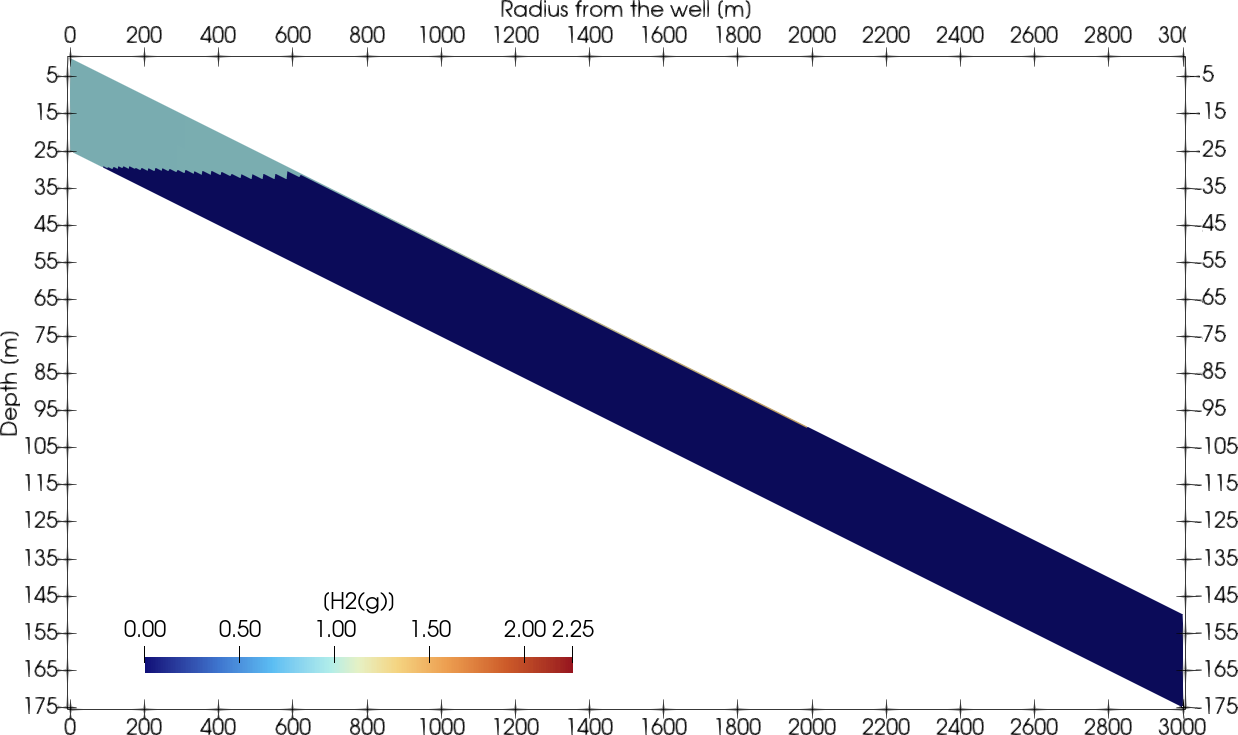}
	}
	\subfloat[molar concentration of $\hdg$ at $t=1000$ days]{
		\includegraphics[width=\2dimgw]{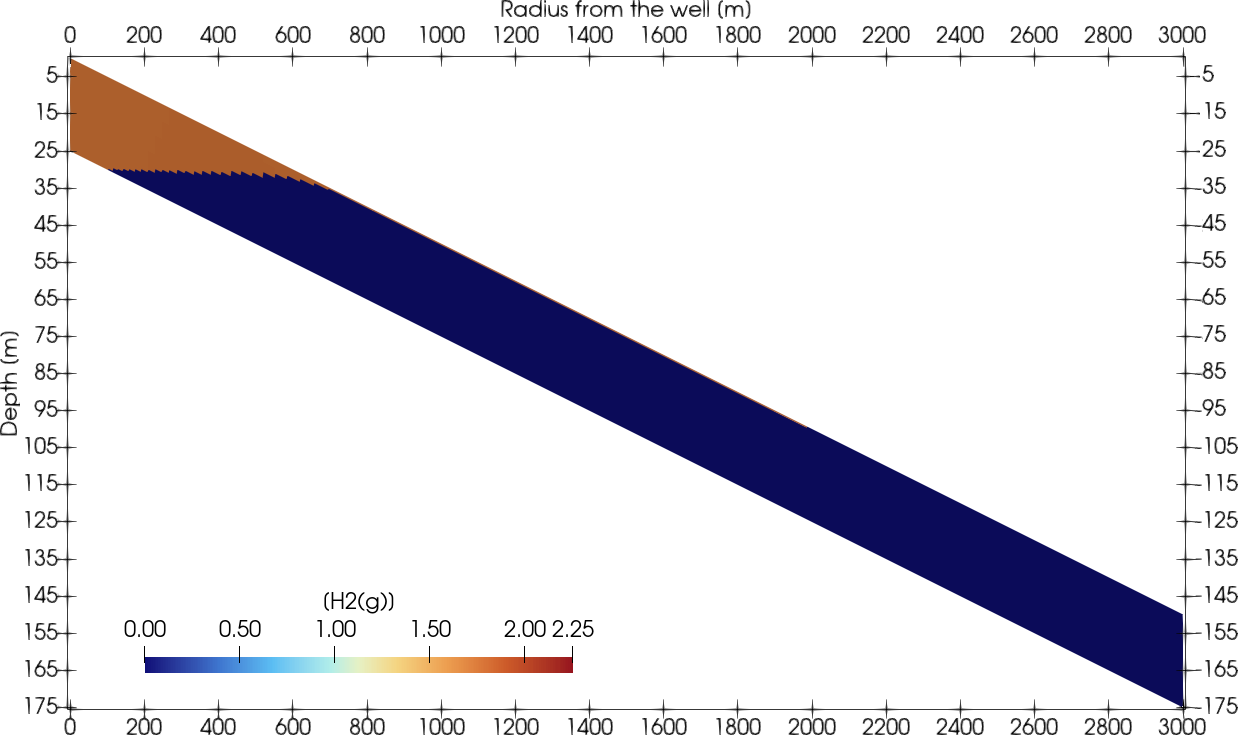}
	}
	\quad
	\subfloat[molar concentration of $\hdsg$ at $t=900$ days]{
		\includegraphics[width=\2dimgw]{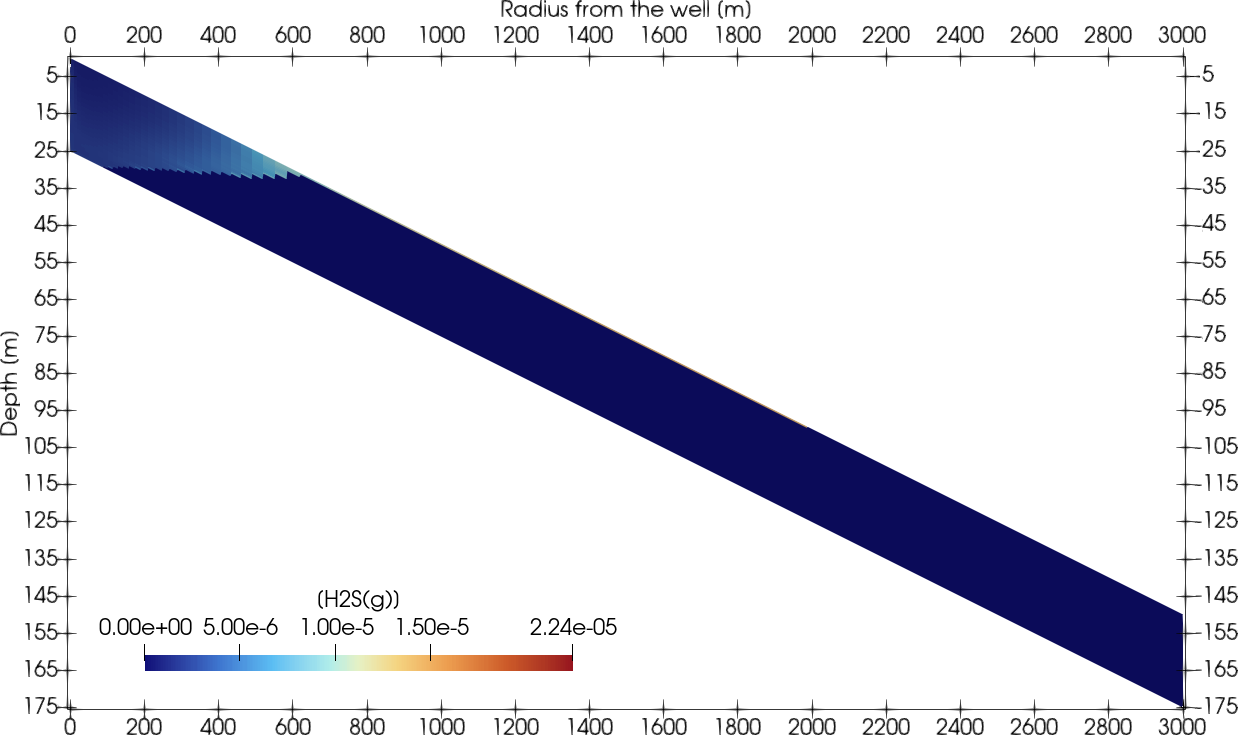}
	}
	\subfloat[molar concentration of $\hdsg$ at $t=1000$ days]{
		\includegraphics[width=\2dimgw]{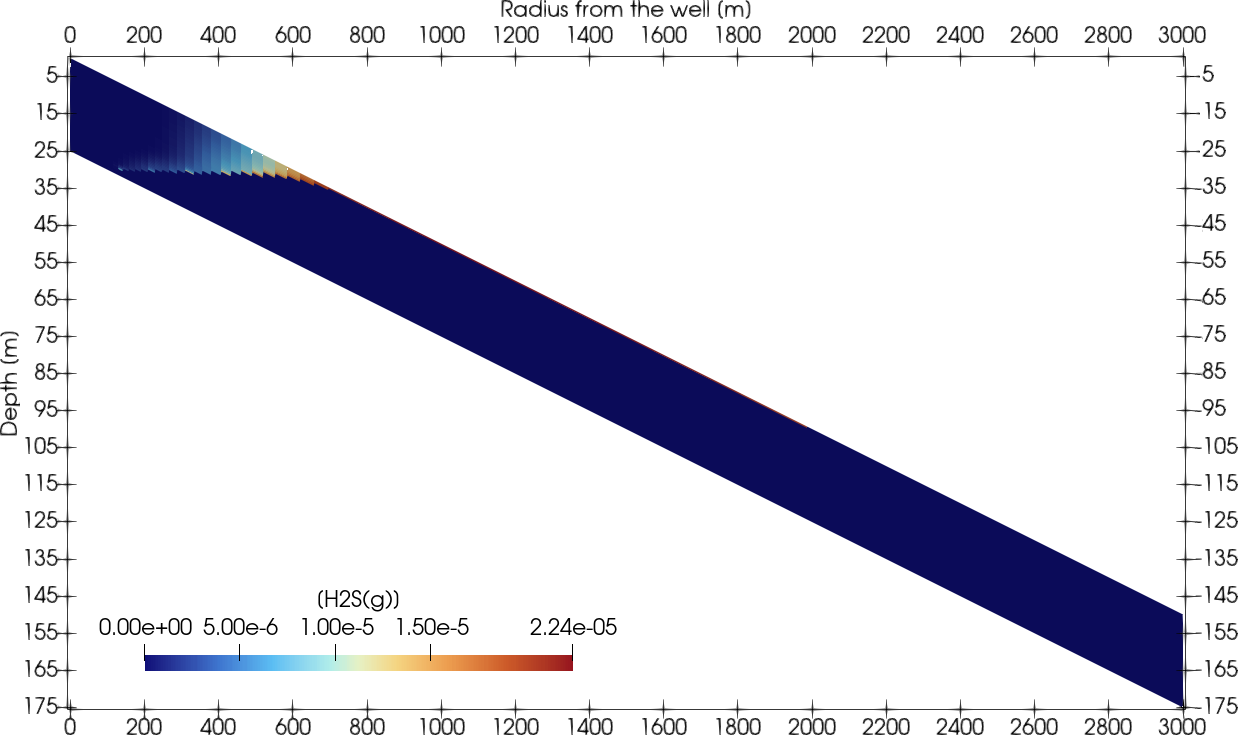}
	}
	\quad
	\subfloat[pH at $t=900$ days]{
		\includegraphics[width=\2dimgw]{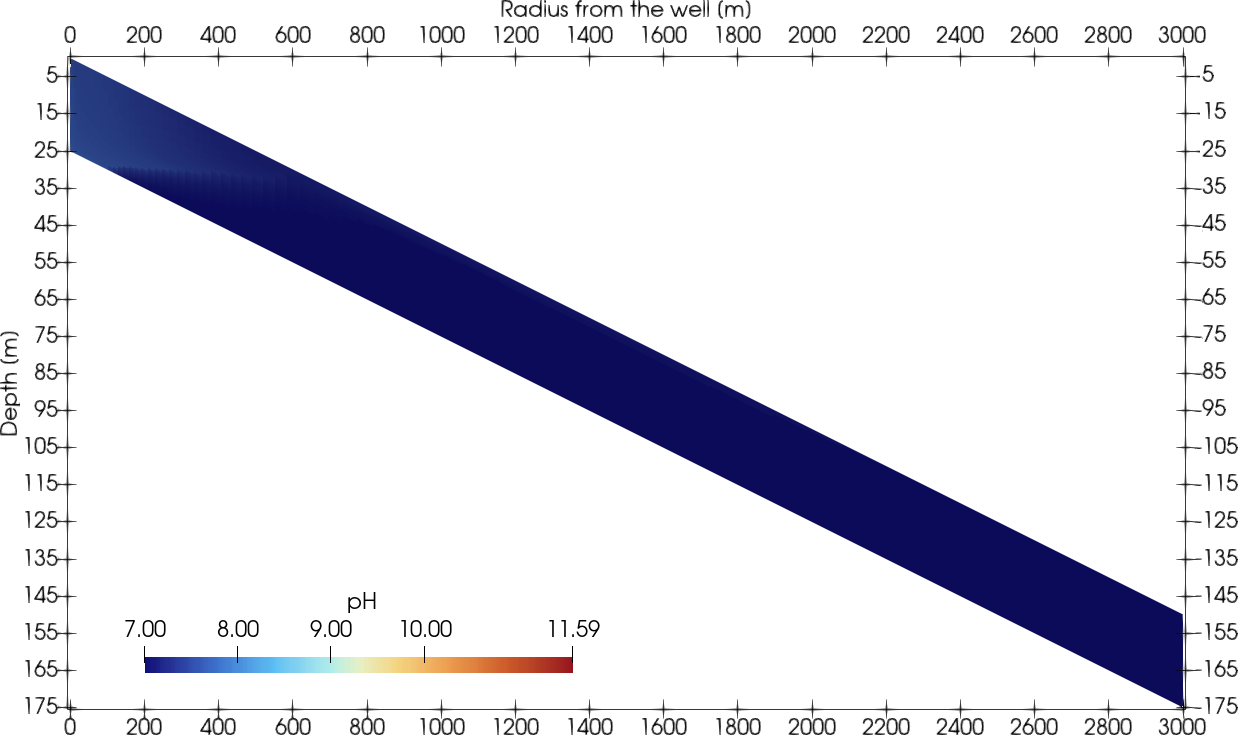}
	}
	\subfloat[pH at $t=1000$ days]{
		\includegraphics[width=\2dimgw]{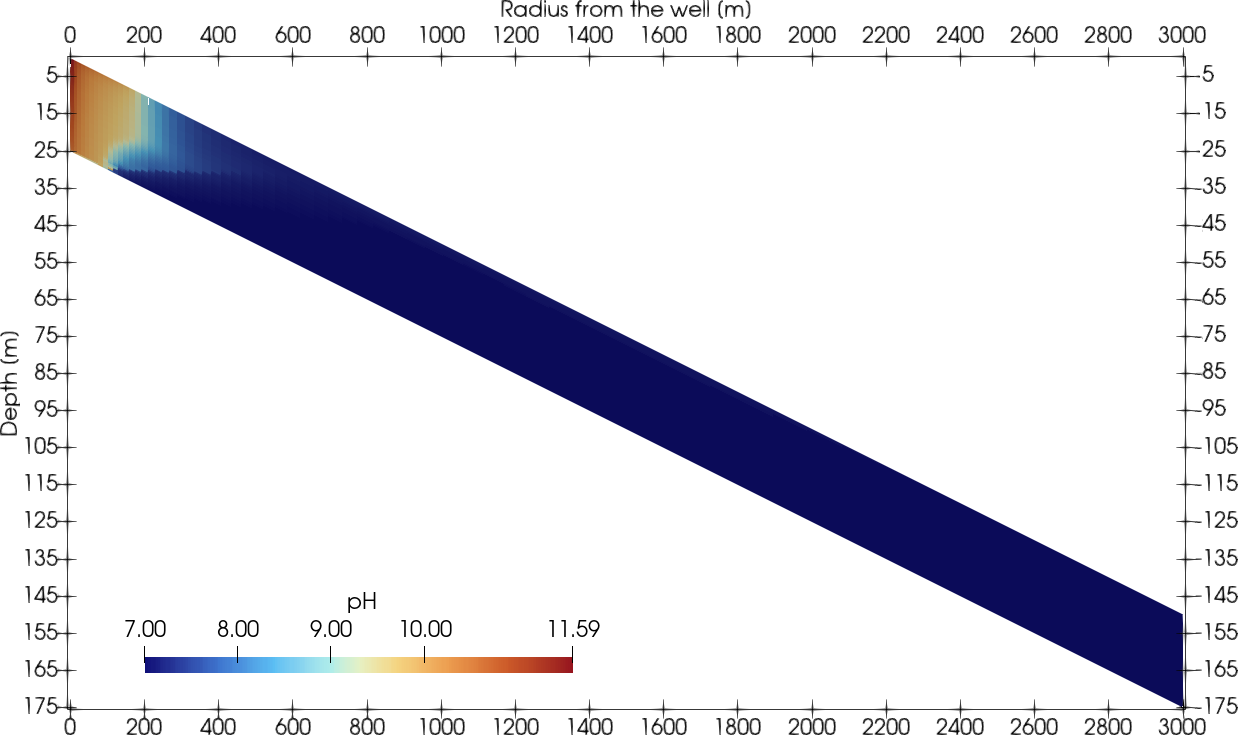}
	}
	\caption{Results of the gaseous $\hd$ injection and withdrawal simulation in a reservoir with radial geometry at times $t = 900$ days and $t = 1000$ days. The molar concentrations of $\hdg$ and $\hdsg$, as well as pH are shown.}
	\label{fig:ben-hd-rad-last}
	\end{figure}

	The simulations are performed on the mesh shown in Figure~\ref{fig:h2-mesh}. Figure~\ref{fig:ben-hd-rad-first} and Figure~\ref{fig:ben-hd-rad-last} show the formation of $\hdsg$ at the front of $\hdg$, as expected. A substantial increase in aqueous pH occurs in the region affected by hydrogen injection. This behavior is likely related to the simplified chemical system considered, which does not account for additional buffering mechanisms, particularly those associated with clay minerals. Although carbonate equilibria and calcite are included, they are insufficient to prevent such a pronounced pH increase in the present simulations.
	
	\section*{Conclusion}
	
	A multiphase reactive transport model based on Morel's formalism and Coats' formulation has been developed, taking into account both equilibrium and kinetic reactions. The system is based on an extension of Coats' formulation to reactive transport combined with a fully implicit global approach. Phase disappearance, including for the aqueous phase, is rigorously handled by the introduction of additional unknowns selected from the restriction to present phases of Morel's formula matrix. These additional unknowns not only avoid the singularity of the system, but they also track  the possible appearance of absent primary species from neighboring cells. The numerical model has been validated and assessed on a recent $\cod$ benchmark in a 2D configuration including gravitational effect and aqueous phase disappearance in the vicinity of the injection well. The use of flash-based phase stability criteria significantly improves the robustness of Newton's method by reducing oscillations associated with phase appearance and disappearance during nonlinear iterations. Finally, it has been successfully tested on a 2D axisymmetric simulation of a hydrogen ($\hd$) injection and production cycles scenario.
	
	Future work will focus on improving the efficiency and applicability of the model by incorporating clay reactions to buffer pH, including relevant minerals among the primary species.
	
	\paragraph{Acknowledgements}
	
	The authors would like to thank STORENGY for funding this work and allowing its publication, and David Segura for fruitful discussions on the $\hd$ storage test case.
	
	\appendix
	
	\section{List of symbols}
	\label{sec:notations}
	
	We refer to Table~\ref{tab:symbols} for a non-exhaustive list of symbols throughout the article.
	
	\begin{table}[H]
	\centering
	\renewcommand{\arraystretch}{1.4}
	\begin{tabular}{|c|c|c|}
		\hline
		\multicolumn{3}{|>{\columncolor{gray!40}}c|}{\textbf{Reservoir}} \\
		\hline
		\rowcolor{gray!20}
		Symbol & Description & Unit \\
		\hline
		$S^\beta$ & Saturation of phase $\beta$ & [$-$] \\
		$\eta^\beta$ & Mobility of phase $\beta$ & [$\mathrm{Pa^{-1} \cdot s^{-1}}$] \\
		$\zeta^\beta$ & Volumetric molar density of phase $\beta$ & [$\mathrm{mol \cdot m^{-3}}$] \\
		$\phi$ & Volume fraction of fluid phases (porosity) & [$-$] \\
		$C_e$ & Mole fraction of chemical species $e$ & [$-$] \\
		$n_e$ & Volumetric molar concentration of chemical species $e$ & [$\mathrm{mol \cdot m^{-3}}$] \\
		$\phi^s$ & Volume fraction of reactive rock & [$-$] \\
		$1-\phi_T$ & Volume fraction of non-reactive rock & [$-$] \\
		$\mathbb{K}$ & Permeability tensor & [$\mathrm{m^2}$] \\
		$k_\text{r}^\beta$ & Relative permeability of phase $\beta$ & [$-$] \\
		$\mu^\beta$ & Viscosity of phase $\beta$ & [$\mathrm{Pa \cdot s}$] \\
		$\u^\beta$ & Velocity of phase $\beta$ & [$\mathrm{m \cdot s^{-1}}$] \\
		$P^\beta$ & Pressure of phase $\beta$ & [$\mathrm{Pa}$] \\
		$\rho^\beta$ & Volumetric mass density of phase $\beta$ & [$\mathrm{kg \cdot m^{-3}}$] \\
		\hline\hline
		\multicolumn{3}{|>{\columncolor{gray!40}}c|}{\textbf{Well}} \\
		\hline
		\rowcolor{gray!20}
		Symbol & Description & Unit \\
		\hline
		$P^\mathrm{w}_\mathrm{bh}$ & Bottom-hole pressure & [$\mathrm{Pa}$] \\
		$Q^{\mathrm{r}\to\mathrm{w}}_h$ & Discretized total molar flow rate & [$\mathrm{mol \cdot s^{-1}}$] \\
		$\overline{P}$ & Target bottom-hole pressure & [$\mathrm{Pa}$] \\
		$\overline{Q}$ & Target discretized total molar flow rate & [$\mathrm{mol \cdot s^{-1}}$] \\
		$\mathrm{WI}$ & Well index & [$\mathrm{m}$] \\
		\hline
	\end{tabular}\vspace{0.2cm}	
	\caption{Non-exhaustive list of symbols used through the article.}
	\label{tab:symbols}
	\end{table}	
	
	\section{Parameters of the $\hd$ test case}
	\label{sec:values-h2}
	
	Table~\ref{tab:param-h2}, Table~\ref{tab:hd-init}, Table~\ref{tab:hd-const} and Table~\ref{tab:mol-mass-h2} give respectively the parameters values, the initial conditions of the reservoir, the reaction equilibrium constants and the molar masses of the chemical species used for the $\hd$ storage test case.
	
	\begin{table}[H]
	\centering
	\renewcommand{\arraystretch}{1.4}
	\begin{tabular}{|c|c|c|}
		\hline
		\rowcolor{gray!40}
		Symbol & Value & Unit \\
		\hline
		$\phi_0$ & $0.27$ & [$-$] \\
		$\mathbb{K}_0$ & $k_0\,\mathrm{Id}$, with $k_0=5\times 10^{-12}$ & [$\mathrm{m^2}$] \\
		$\mathbb{K}$ & Kozeny-Carman \eqref{eq:kozeny-carman} & [$\mathrm{m^2}$] \\
		$k_\text{r}^\l$, $k_\text{r}^g$ & Brooks-Corey \eqref{eq:brooks-corey} & [$-$] \\
		$\mu^\l$ & $10^{-3}$ & [$\mathrm{Pa \cdot s}$] \\
		$\mu^g$ & $8.9\times 10^{-6}$ & [$\mathrm{Pa \cdot s}$] \\
		$P_0$ & $31\times 10^5$ & [$\mathrm{Pa}$] \\
		$r_\text{w}$ & $0.05$ & [$\mathrm{m}$] \\
		$\rho^\l$ & Linear compressibility \eqref{eq:lin-compress-fluid} & [$\mathrm{kg\cdot m^{-3}}$] \\
		$\rho^g$ & Ideal gaz mixture \eqref{eq:ideal-gaz} & [$\mathrm{kg\cdot m^{-3}}$] \\
		$\zeta^\l$, $\zeta^g$ & Equation~\eqref{eq:zeta-rho} with Table~\ref{tab:mol-mass-h2} & [$\mathrm{mol\cdot m^{-3}}$] \\
		\hline
	\end{tabular}\vspace{0.2cm}	
	\caption{$\hd$ storage test case: Value of the parameters of the model.}
	\label{tab:param-h2}
	\end{table}	
	
	\begin{table}[H]
	\centering
	\begin{tabular}{|c|c|c|c|c|}
		\hline
		\rowcolor{gray!40} $\hdl$ & $\soqdm$ & $\hdsl$ & DIC & pH \\
		\hline
		$0$ & $20$ $\mathrm{mg \cdot L^{-1}}$ & $0$ & $5\times 10^{-4}$ $\mathrm{mol \cdot L^{-1}}$ & $7$ \\
		\hline
	\end{tabular}\vspace{0.2cm}		
	\caption{$\hd$ storage test case: The initial state of the reservoir is assumed to be at chemical equilibrium, where DIC stands for \textit{Dissolved Inorganic Carbon}, defined as $\text{DIC} = \codl + \hcotm + \cotdm$.}
	\label{tab:hd-init}
	\end{table}
	
	\begin{table}[H]
	\centering
	\begin{tabular}{|c|c|c|c|c|c|c|c|}
		\hline
		\rowcolor{gray!40} $K_{\hsm}$ & $K_{\ohm}$ & $K_{\hcotm}$ & $K_{\cotdm}$ & $K_{\hdsg}$ & $K_{\hdg}$ & $K_{\codg}$ & $K_{\cacot}$ \\
		\hline
		$-6.989$ & $-13.997$ & $-6.354$ & $-16.683$ & $5.785$ & $753.908$ & $17.391$ & $8.203$ \\
		\hline
	\end{tabular}\vspace{0.2cm}		
	\caption{$\hd$ storage test case: Constants associated with the equilibrium chemical reactions ($\log_{10}$ scale).}
	\label{tab:hd-const}
	\end{table}
	
	\begin{table}[H]
	\centering
	\begin{tabular}{|c|c|}
		\hline
		\rowcolor{gray!40} Molar mass & Value [$\mathrm{kg\cdot mol^{-1}}$] \\
		\hline
		$M_{\hdsl}$, $M_{\hdsg}$ & $34.0809\times 10^{-3}$ \\
		$M_{\hdl}$, $M_{\hdg}$ & $2.01588\times 10^{-3}$ \\
		$M_{\codl}$, $M_{\codg}$ & $44.0095\times 10^{-3}$ \\
		$M_{\hdo}$ & $18.01528\times 10^{-3}$ \\
		$M_{\hp}$ & $1.00794\times 10^{-3}$ \\
		$M_{\cadp}$ & $40.078\times 10^{-3}$ \\
		$M_{\nap}$ & $22.9892207\times 10^{-3}$ \\
		$M_{\soqdm}$ & $96.0637\times 10^{-3}$ \\
		$M_{\hsm}$ & $33.0729\times 10^{-3}$ \\
		$M_{\ohm}$ & $17.00734\times 10^{-3}$ \\
		$M_{\hcotm}$ & $61.0174\times 10^{-3}$ \\
		$M_{\cotdm}$ & $60.01\times 10^{-3}$ \\
		\hline
	\end{tabular}\vspace{0.2cm}		
	\caption{$\hd$ storage test case: Molar masses of the chemical species.}
	\label{tab:mol-mass-h2}
	\end{table}
	
	\printbibliography
	
\end{document}